\documentclass{amsart}
\usepackage{epsfig,amsmath,latexsym}
\usepackage{amsfonts,amssymb,graphicx}
\usepackage{amscd, amsthm, amsmath,eucal, verbatim}

\newcommand{\bproof}{\paragraph {\bf Proof.}}

\newcommand{\bv}{{\bf v}}
\newcommand{\bw}{{\bf w}}
\newcommand{\RR}{{\mathbb{R}}}
\newcommand{\ZZ}{{\mathbb{Z}}}
\newcommand{\sC}{{\mathcal{C}}}
\newcommand{\sH}{{\mathbb{H}}}

\newtheorem{theorem}{Theorem}
\newtheorem{lemma}[theorem]{Lemma}
\newtheorem{corollary}[theorem]{Corollary}
\newtheorem{claim}{Claim}

\begin{document}
\title{The Computational Complexity of Knot Genus and Spanning Area}
\date{\today}
\author{Ian Agol}
\address{Dept. of Mathematics, Statistics, and Computer Science, 
University of Illinois at Chicago, 
Chicago, IL 60607}
\email{agol@math.uic.edu}
\thanks{Partially supported by ARC grant 420998.}
\author{Joel Hass}
\address{School of Mathematics, Institute for Advanced Study and
Department of Mathematics, University of California, Davis, California 95616}
\email{hass@math.ucdavis.edu}
\thanks{This work was carried out while the second author was visiting the
Institute for Advanced Study, partially supported by NSF grant DMS-0072348,
and by a grant to the Institute for Advanced Study by AMIAS.}
\author{William Thurston}
\address{Department of Mathematics, University of California, Davis
California 95616}
\thanks{Third author partially supported by NSF grant DMS-9704286.}
\email{wpt@math.ucdavis.edu}
\begin{abstract}
We show that the problem of
deciding whether a polygonal knot in a closed
three-dimensional manifold bounds a surface of genus at most $g$,
is {\bf NP}-complete.
We also show that the problem of deciding whether a curve in a PL manifold
bounds a surface of area less than a given constant $C$ is
{\bf NP}-hard.
\end{abstract}
\maketitle

\section{Introduction} \label{intro}
In this paper we investigate the computational complexity
of some problems in three-dimensional topology and geometry.
We show that the problem of
determining a bound on the genus
of a knot in a 3-manifold, is {\bf NP}-complete.
Using similar ideas, 
we show that deciding
whether a curve in a metrized PL 3-manifold
bounds a surface of area less than a given constant $C$
is {\bf NP}-hard.

Determining whether a given knot is trivial or not is
one of the historically central questions in topology.
The problem of finding an algorithm to determine
knot triviality was posed by Dehn \cite{Dehn:10}.
Dehn's investigations into this area led to the formulation
of the word and isomorphism problems, which played an
important role in the development of the theory of
algorithms.
The first algorithm for the unknotting problem was given by Haken \cite{Hak:61}. Haken's procedure is
based on normal surface theory, a method of representing surfaces
introduced by Kneser \cite{Knes:29}.
Analysis of the computational complexity
of this algorithm is more recent. Hass, Lagarias and
Pippenger showed that
Haken's unknotting algorithm runs in time at most $c^t$,
where the knot $K$ is embedded in the 1-skeleton of a
triangulated manifold $M$ with $t$ tetrahedra,
and $c$ is a constant independent of $M$ or $K$ \cite{HLP:99}.
It was also shown in \cite{HLP:99} that the unknotting problem
is in {\bf NP}.

The notion of genus was defined by Seifert \cite{Sei:35}
in 1935 for knots in the 3-sphere, and
extends directly to knots in an arbitrary 3-manifold $M$.
Given a knot $K$, consider the class ${\mathcal S}(K)$ of all orientable
spanning surfaces for $K$. These are surfaces embedded in $M$ with
a single boundary component that coincides with $K$.
Seifert showed that this class is non-empty for any
knot $K$ in the 3-sphere.  For knots in a general manifold,
${\mathcal S}(K)$ is non-empty when $K$ represents a trivial element
in the first integer homology group of $M$.
The genus $g(K)$ of a knot $K$ is
the minimum genus of a surface in ${\mathcal S}(K)$, or $\infty$ if
${\mathcal S}(K) = \emptyset$.
The genus measures one aspect of the degree of ``knottedness'' of a curve.

The unknotting problem is a special case of the more general
problem of determining the genus of a knot in a 3-manifold. Given a knot $K$ in
an orientable 3-manifold and a positive integer $g$, this problem asks
for a procedure to determine whether the {\em knot genus} of $K$, the minimal genus
of an orientable spanning surface for $K$ in a 3-dimensional manifold, is at most $g$.
A knot is trivial, or unknotted, precisely when its genus is zero.
We will show that the problem of
determining the genus of a knot in a 3-manifold is {\bf NP}-complete.
Previous results on this problem were given
in \cite{HLP:99}, where it was shown to lie in {\bf PSPACE}, roughly
the class of problems that run in polynomial space.
No lower bounds on the running time were previously known.

We work with 3-manifolds that are triangulated and orientable, and
with orientable embedded surfaces. This is not a significant
restriction, since all
compact 3-dimensional manifolds admit unique PL structures \cite{Moise:52}.
A {\em knot} in a triangulated 3-manifold $M$ is
a connected simple (non self-intersecting) closed
curve in the 1-skeleton of $M$.
Any smooth knot in a smooth manifold, or more generally any tame knot, is equivalent
to a knot that lies in the 1-skeleton of some triangulation.

We formulate the problem of computing
the genus as a language-recognition problem in the usual way,
see \cite{GJ:79}.
In 1961 Schubert \cite{Sch:61}, in an extension of Haken's work,
showed the decidability of the problem:

\begin{tabbing}
Problem: \= 3-MANIFOLD KNOT GENUS \\
INSTANCE: \=A triangulated 3-dimensional manifold $M$, a knot
$K$ in the\\
1-skeleton of $M$, and a natural number $g$. \\
QUESTION: \=Does the knot $K$ have $g(K) \le g$?
\end{tabbing}

The size of an instance is measured by the sum of the
number of tetrahedra in $M$.
In Section~\ref{sectionnph} we establish

\begin{theorem} \label{thmnph}
{\sc 3-MANIFOLD KNOT GENUS} is {\bf NP}-hard.
\end{theorem}

It was established in \cite{HLP:99} that {\sc 3-MANIFOLD KNOT GENUS} is in
{\bf PSPACE}.  We
improve this bound in Section~\ref{sectiongnp} .

\begin{theorem} \label{thmgnp}
{\sc 3-MANIFOLD KNOT GENUS} is {\bf NP}.
\end{theorem}

In combination these two results give:

\begin{theorem} \label{3mkgnpc}
{\sc 3-MANIFOLD KNOT GENUS} is {\bf NP}-complete.
\end{theorem}

Theorem~\ref{thmnph} is proved through a connection to
{\sc ONE-IN-THREE SAT}, a known  {\bf NP}-complete problem that will be reviewed
in Section~3.
The theorem carries out a construction that
transforms an instance of {\sc ONE-IN-THREE SAT} to an instance of {\sc 3-MANIFOLD KNOT GENUS}.
By ``transform'' we mean that an instance of one problem is changed
to an instance of the second by a procedure that requires time polynomial
in the size of the instance.
To a boolean expression representing an instance of  {\sc ONE-IN-THREE SAT}
we associate a positive integer $g$
and a certain knot in a triangulated, compact 3-manifold.
This knot bounds a surface of genus at most $g$ exactly when there
is a truth assignment to the boolean expression satisfying the
requirements of  {\sc ONE-IN-THREE SAT}.
Since {\sc ONE-IN-THREE SAT} is {\bf NP}-hard,
this establishes that {\sc 3-MANIFOLD KNOT GENUS}
is also {\bf NP}-hard.

In Section~\ref{sectiongnp} we prove Theorem~\ref{thmgnp}, giving a
certificate which demonstrates in polynomial time that a genus $g$ knot $K$
bounds a surface of genus at most $g$.
The argument in \cite{HLP:99}
established that the unknotting problem is {\bf NP}
using the existence of a normal disk that
lies along an extremal ray in the space of normal solutions,
called a vertex surface in Jaco-Tollefson \cite{JT:95}.
The existence of such an extremal
normal surface of minimal genus spanning a knot
is not known, so a new technique is needed (See \cite{JO:84} for known results here). 
This is provided in Theorem~\ref{orbitcount}, which gives
an algorithm to count the number of orbits of a type of pseudogroup
action on a set.
Theorem~\ref{orbitcount} seems likely to have more general applicability.
In Section~\ref{pseudogroup} we describe this algorithm and in
Section~\ref{sectiongnp} we
apply it to a pseudogroup action that arises in the theory of normal
surfaces.  This allows us to determine in polynomial time the number
of components in a normal surface described by an integer vector
in ${\ZZ_+^{7t}}$. In particular we are able to certify that a normal
surface is connected, orientable and has connected boundary.
Since calculating the Euler characteristic of a normal
surface can be done efficiently,
establishing orientability and
connectedness are the key steps in constructing 
a certificate of its genus.

In Section~\ref{extended} we extend the orbit counting algorithm to
allow the counting of additional integer weight sums associated to
each orbit.
This allows for the polynomial time calculation of the genus of 
all the components of a
fundamental normal surface, as well as a count of
the number of components.

The genus and the area of a surface are closely connected.
In Section~\ref{area} we extend the methods developed in
studying genus to study the problem of determining the smallest
area of a spanning surface for a curve in a 3-manifold.
We show that computing an upper bound on the
area of a smallest area spanning surface is  {\bf NP}-hard.

We refer to \cite{Papadimitriou:94} for a
discussion of complexity classes such as {\bf NP} and {\bf PSPACE} and
\cite{{Wel:93a}} for a discussion of complexity problems in
low-dimensional topology.

{\bf Remarks:}
\begin{enumerate}
\item Knots are often studied in $\RR^3$ or $S^3$ rather
than in a general manifold. Our methods show that
determining knot genus in $\RR^3$ or $S^3$, or any
fixed manifold, is {\bf NP}. It is not clear whether
the corresponding problem remains {\bf NP}-hard  
if one restricts consideration to knots in $\RR^3$ or $S^3$.

\item Casson has shown that a procedure to determine whether a 3-manifold
is homeomorphic to the 3-sphere,
following the 3-sphere recognition
algorithm described in \cite{Rubinstein:94}
and \cite{Thompson:94}, runs in time
less than $3^tp(t)$, where $p(t)$ is a polynomial.
In the direction of lower bounds, it was shown in \cite{JVW:90}
that determining certain values of
the Jones polynomial of alternating links
is {\bf \#P}-hard.
\end{enumerate}

We are grateful to the referee for numerous suggestions on 
the exposition of this paper.

\section{Normal Surfaces} \label{ns}

\begin{figure}[hbtp]
\centering
\includegraphics[width=.6\textwidth]{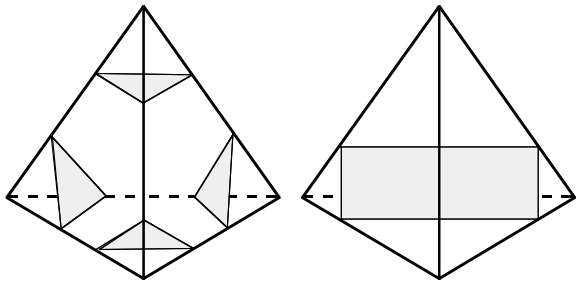}
\caption{
\label{normal}
Elementary disks in a normal surface.}
\end{figure}

General surfaces in 3-manifolds can wind and
twist around the manifold in complicated ways.
Kneser described a procedure in which surfaces can be ``pulled taut'',
until they take a simple and rigid position  \cite{Knes:29}.
In this {\em normal} position,
they have very succinct algebraic descriptions.
We use an approach to normal surfaces in triangulated 3-manifolds based on work of
Jaco-Rubinstein \cite{JR:89} and Jaco-Tollefson
\cite{JT:95}.
A {\em normal surface} $S$ in a triangulated compact 3-manifold
$M$ is a $PL$-surface whose intersection with each tetrahedron in $M$
consists of a finite number of disjoint {\em elementary disks}. These are properly
embedded disks that are isotopic to either
triangles or quadrilaterals as shown in Figure~\ref{normal}, by an isotopy
preserving each face of the tetrahedron.

Within each tetrahedron of $M$ there are four
possible triangles and three possible quadrilaterals, up to a
{\em normal isotopy of $M$}, an isotopy which leaves
each cell of the triangulation of $M$ invariant.
W. Haken observed that
a normal surface is determined up to
such isotopies by the number of
pieces of each of the seven kinds of elementary disks that
occur in each tetrahedron, or a vector in ${\ZZ_+^{7t}}$.
A normal surface $S$ is described by a non-negative integer vector
$\bv = \bv (S) \in \ZZ^{7t}$, that gives the
{\em normal coordinates} of $S$.
There is a homogeneous rational cone
${\sC}_M$ in $\RR^{7t}$, called the {\em Haken normal cone},
that contains the vectors $\bv (S)$ for all normal surfaces $S$ 
in $M$.

If $\bv = ( v_1 , v_2 , \ldots , v_{7t} ) \in \RR^{7t}$,
then the Haken normal cone is specified by
linear equations and inequalities of the form
$$
v_{i_1} + v_{i_2} = v_{i_3} + v_{i_4} \quad \mbox{(up to $6
t$ equations)}~,
$$
$$
v_i \geq 0~, \quad \mbox{for}~~ 1 \leq i \leq 7t~.
$$
The first set of equations expresses
{\em matching conditions,} which say that the number of edges
on a common triangular face of two adjacent tetrahedra,
coming from a collection of elementary disks in each of
the tetrahedra, must match. For each triangular face there
are three types of edges (specified by a pair of edges on the
triangle), which yield 3~matching conditions per face. 
Triangular faces in the boundary $\partial M$ 
give no matching equations.
The second set are called the {\em positivity} conditions.
The cone ${\sC}_M$ is {\em rational,}
because the above equations
have integer coefficients. We let
$ {\sC}_M (\ZZ ) = {\sC}_M \cap \ZZ^{7t}$
denote the set of integral vectors in the cone
${\sC}_M$. An additional set of conditions, the
{\em quadrilateral conditions,}
is required for an integral vector in the cone to
correspond to the normal coordinates of an embedded surface.
This condition states that of the three types
of quadrilateral found in each tetrahedron,
only one can occur in the vector with non-zero coefficient.
The quadrilateral conditions are
required because two distinct types of quadrilateral in
a single tetrahedron necessarily intersect, and we are interested
in embedded, non-self-intersecting surfaces.
A vector in the Haken normal cone
that satisfies the quadrilateral conditions
corresponds to an embedded normal surface.
This surface is unique up to a
normal isotopy.
However it is important to note that the surface corresponding
to a given vector in the cone $ {\sC}_M$ may not be
connected.
A normal surface meets each edge of a triangulation in a finite number
of points.  The sum of these intersection numbers over all the
edges of a triangulation is called the {\em weight} of the normal surface.

A {\em fundamental normal surface} is a normal surface $S$
such that
\begin{displaymath}
\bv (S) \ne \bv_1 + \bv_2 ~, \quad \mbox{with $\bv_1 , \bv_2 \in
{\sC}_M ( \ZZ )
\backslash \{ \bf 0 \} $ }.
\end{displaymath}
In the terminology of integer programming,
such a vector $\bv (S)$ is an element of the
{\em minimal Hilbert basis}
${\sH}({\sC}_M )$ of ${\sC}_M$,
see Schrijver \cite[Theorem~16.4]{Schrijver:86}.
A fundamental normal surface is always
connected, but connected normal surfaces need not be fundamental.
A {\em  vertex minimal solution} is a special kind of fundamental surface, one that
corresponds to a solution 
of the normal surface equations that lies along an extremal ray of the
cone of solutions and is not a multiple of another such extremal solution.
Hass, Lagarias and Pippenger \cite[Lemma 6.1]{HLP:99} gave a bound for
the size of the vectors corresponding to any fundamental surface.

\begin{theorem} \label{bound}
Let $M$ be a triangulated compact 3-manifold, possibly with boundary,
that contains $t$ tetrahedra.
\begin{itemize}
\item
Any vertex minimal solution $\bv \in \ZZ^{7t}$ of the Haken
normal cone $C_M$ in $\RR^{7t}$ has
$\displaystyle \max_{1 \le i \le 7t} (v_i) \le 2^{7t-1} ~.
$
\item
Any minimal Hilbert basis element $\bv \in \ZZ^{7t}$ of the
 Haken
fundamental cone $C_M$ has
$\displaystyle \max_{1 \le i \le 7t} (v_i) < t \cdot 2^{7t+
2} ~. $
\end{itemize}
\end{theorem}

Schubert \cite{Sch:61} showed that a surface of smallest genus
spanning $K$ can be found among the fundamental surfaces.

\begin{theorem} \label{mingenusfund}
There is a minimal genus spanning surface for $K$ which is a
fundamental normal surface.
\end{theorem}

Similar to the theory of normal surfaces, though somewhat
easier, is the theory of {normal curves}.
These are curves on a surface that intersect
each triangle in a collection of {\em normal arcs},
arcs that have endpoints on distinct edges. 
Normal curves arise as the boundaries
of normal surfaces in a manifold with boundary.
Since there are three such arcs in each triangle, 
normal isotopy classes of normal curves in a surface that contains $t$ triangles
are described by integer vectors in ${\ZZ_+^{3t}}$.

\section{{\sc 3-MANIFOLD KNOT GENUS} is {\bf NP}-hard} \label{sectionnph}

In this section we show how to
reduce an instance of {\sc ONE-IN-THREE SAT}
to an instance of {\sc 3-MANIFOLD KNOT GENUS}.
Since {\sc ONE-IN-THREE SAT} is known to be {\bf NP}-hard,
this establishes that {\sc 3-MANIFOLD KNOT GENUS} is also {\bf NP}-hard.
The problem {\sc ONE-IN-THREE SAT} concerns
logical expressions involving collections of literals (boolean variables or their negations)
gathered in clauses consisting of three literals connected with $\vee$'s.  The logical expression
contains a collection of clauses connected with $\wedge$'s.

\begin{tabbing}
Problem:  ONE-IN-THREE SAT\\
INSTANCE: A set $U$ of variables and a collection $C$ of clauses over $U$\\
such that each clause $c \in C$ contains 3 literals. \\
QUESTION: \=Is there a truth assignment for $U$ such that each clause in $C$\\
has exactly one true literal?
\end{tabbing}

Schaefer \cite{Sch:78} established that {\sc ONE-IN-THREE SAT} is
{\bf NP}-complete.
To prove Theorem~\ref{thmnph}, establishing that
{\sc 3-MANIFOLD KNOT GENUS} is {\bf NP}-hard, we
show that an arbitrary problem in
{\sc ONE-IN-THREE SAT} can be reduced
in polynomial time to a problem in {\sc 3-MANIFOLD KNOT GENUS}.
See Garey and Johnson \cite{GJ:79} for a discussion and many examples of such reductions.

\begin{figure}[hbtp]
\centering
\includegraphics[width=.4\textwidth]{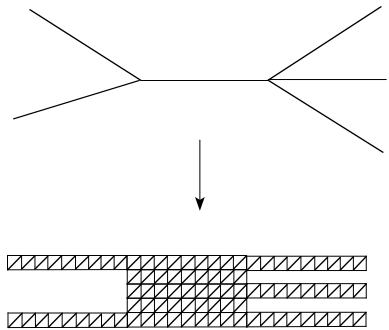}
\caption{A branching surface $B$, shown in cross-section,
is ``thickened'' to produce a triangulated 3-manifold with boundary.}
\label{thicken}
\end{figure}

Let $U = \{ u_1, u_2, \dots u_n \}$ be a set of variables and
$C = \{ c_1, c_2, \dots c_m \}$
be a set of clauses in an arbitrary instance of {\sc ONE-IN-THREE SAT}.
We will describe a knot $K$ in a compact 3-dimensional manifold
(with no boundary) and an integer $g$
such that $K$ bounds a surface of genus smaller or equal to $g$
if and only if $C$ is satisfiable so that each
clause in $C$ contains exactly one true literal.
We construct the 3-manifold $M$ in stages.  First we construct a
2-dimensional simplicial complex, then we thicken this complex,
replacing triangles with subdivided triangular prisms,
getting a triangulated,
3-dimensional manifold with boundary, as indicated in
Figure~\ref{thicken}. Finally we use
a doubling construction, taking two copies of the manifold
with boundary and gluing their boundaries
together, to obtain a closed 3-manifold.

\begin{figure}[hbtp]
\centering
\includegraphics[width=.5\textwidth]{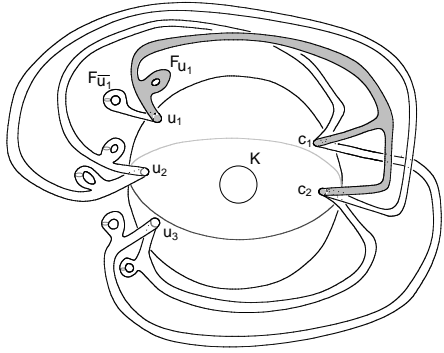}
\caption{A branching surface with boundary curve $K$
corresponding to the boolean expression
$(u_1 \vee u_2 \vee u_3) \wedge (u_1 \vee \bar u_2 \vee \bar u_3)$.
The shaded surface $F_{u_1}$ indicates
the occurrence of $u_1$ in each of the clauses $c_1$ and $c_2$.}
\label{branching}
\end{figure}

To begin, we form a type of singular surface $B$, that we call
a {\em branching} surface, by identifying
boundary curves of a collection of $2n+1$ surfaces with boundary, each
forming what we refer to as a {\em piece} of the branching surface.
We construct this collection of surfaces as follows.
Let $k_i$ be the number of times that the variable
$u_i$ appears in the collection of clauses,
and ${\bar k}_i$ be the number of times that the negation ${\bar u}_i$ of $u_i$ appears.
For $i=1,\dots,n$, let $F_{u_i}$
and $F_{{\bar u}_i}$ be genus one surfaces with
$k_i + 1$ and ${\bar k}_i + 1$ boundary curves respectively.
Also set $F_0$ to be a planar surface with $n+m+1$ boundary curves.
One of these boundary curves will later become the knot $K$.
The branching surface $B$ is constructed by identifying these surfaces
along appropriate boundary components as indicated in Figure~\ref{branching}.

Branching occurs when more than
two boundary curves are identified along a single curve.
We identify pairs of boundary curves by giving a homeomorphism between them.
Up to isotopy, this is determined by specifying an
orientation on the curves and setting the homeomorphism to
be orientation reversing. 
We first fix an orientation on each of $F_0, F_{u_i}$
and $F_{{\bar u}_i}$. This induces an orientation on each boundary
curve. All identifications will involve gluing a
boundary component of $F_{u_i}$
or $F_{{\bar u}_i}$ to a boundary component of $F_0$, and we require this
gluing to be an orientation reversing homeomorphism.
Label by $K, u_1, \dots u_{n}, c_1, \dots c_m$
the $1+n+m$ boundary components of $F_0$.
The boundary component $K$ of $F_0$, which will become our
knot, has nothing identified to it.
For each $i$, $1 \le i \le n$, one boundary curve from the
surface $F_{u_i}$ is identified to $u_i$.
The remaining $k_i$ boundary components of $F_{u_i}$
are identified with $k_i$ of the curves $c_{1}, \dots c_{m}
$ on $\partial F_0$,
with one component of $\partial F_{u_i}$ identified with $c_j$
for each occurrence of the literal $u_i$
in the $j^{th}$ clause in $C$.
Similarly, one curve of $\partial F_{\bar{u}_i}$ is
identified to the component of $F_0$ labeled $u_i$, and
the remaining ${\bar k}_i$ boundary components
are glued to  $c_{1}, \dots c_{m}$,
with a component glued to $c_j$ for each
occurrence of the literal ${\bar u}_i$ in the $j^{th}$ clause of $C$.
A total of three surface boundaries are identified
along each of $u_1, \dots , u_n$, and
exactly four surface boundaries are identified
along each of $c_{1}, \dots ,c_{m}$, as in
Figure~\ref{branching}.

\begin{lemma} \label{branch}
There is a truth assignment for $U$ such that each 
clause in $C$ has exactly one true literal if and only if
there is a surface $S$ with connected boundary and genus at most $m+n$
and a continuous map $f: S \to B$ such that $f|_{\partial S}$ is a homeomorphism onto $K$.
\end{lemma}
\bproof
Suppose there is a truth assignment for $U$
such that each clause in $C$ contains exactly one true literal.
Form a surface $S$ inside $B$ by taking the union of $F_0$
and either $F_{u_i}$ if $u_i$ is true, or $F_{{\bar u}_i}$
if $ u_i$  is false.
Then exactly two boundary components will be identified along each of
the boundary components of $F_0$ other than $K$ itself,
and $K$ becomes the boundary of the
resulting embedded surface $S$.
There is a contribution of one to the
genus from each of the literals,
since $F_{u_i}$ and $F_{{\bar u}_i}$ each have genus one,
and a contribution of one to the genus from each handle formed
when a boundary component of a surface
$F_{u_i}$ or $F_{{\bar u}_i}$
is glued to $F_0$ along a curve $c_j$.
The genus of $S$ is therefore equal to $m+n$.

Now suppose there is a surface $S_1$ of genus
$\le n+m$ mapped continuously into $B$ that has a single
boundary component mapped homeomorphically to the boundary curve $K$.
We will show in this case that there is a surface $S_4$ with the
same boundary, consisting of certain pieces of $B$ identified along their
boundaries, having
genus precisely $n+m$, and containing, along each curve $c_j$,
exactly one of the three pieces of surface joining $F_0$.

The map $f_1$ of $S_1$ into $B$ may be quite complicated, winding back and
forth across $B$, but by standard transversality arguments
we can homotop $f_1$ so that it is 
a union of homeomorphisms of subsurfaces of $S_1$ mapped homeomorphically to one of
the pieces $F_0$, $F_{u_i}$, $F_{\bar{u}_i}$ forming $B$.
More precisely, we can perturb $S_1$ by a small
homotopy so that its intersection with the boundary components of $F_0$
is transverse, and pulls back to a collection of simple closed curves
on $S_1$. If any of these curves bounds a disk in $S_1$,
then the disk is mapped into some subsurface $X$ of $B$ while the
boundary curve is mapped to $\partial X$.  Since $X$ is not itself a disk, the
disk can be homotoped into $\partial X$, and
we can therefore homotop $S_1$ in a neighborhood of this disk to remove
a component of $S_1 \cap f^{-1}(\partial F_0)$.
After repeating finitely many times,
each component of the complement of  $S_1 \cap  f^{-1}(\partial F_0)$
in $S_1$ has non-positive Euler characteristic.
The image of $S_1$ in each piece $F_0$, $F_{u_i}$,
$F_{\bar{u}_i}$ has an algebraic degree, which is either even or odd.
This degree equals the number of pre-images in $S_1$
of a generic point in the piece.
The degree of the map from $S_1$ on $F_0$ is odd, since $K$ is the
boundary of $S_1$ and therefore $S_1$ maps an odd number of times to 
points near $K$.
The sum of the degrees 
along each of the pieces meeting a curve $u_i$ or $c_j$ is even,
since $S_1$ has no boundary along these curves.
In particular, for each $1 \le i\le n$, exactly one of $F_{u_i}$,
$F_{\bar{u}_i}$ has odd multiplicity in $S_1$.
Form a new surface $S_2$ by taking the
union of $F_0$ and each of the pieces
$F_{u_i}$, $F_{\bar{u}_i}$ which have odd
multiplicity in $S_1$.
We will show that $\chi (S_2) > \chi (S_1)$.
The surface $S_2$ is obtained from $S_1$ by a series of operations that 
either discard a subsurface with non-positive Euler characteristic or
replace a subsurface of $S_1$ that maps with odd degree to
some piece $F_{u_j}$ with a subsurface mapping homeomorphically to
$F_{u_j}$.
The Euler characteristic of a
discarded subsurface of $S_1$ is smaller or equal to that of the subsurface
that replaces it, so $\chi(S_2) \ge \chi(S_1)$ in either case.

The collection of pieces among $\{ F_{u_i}$, $F_{\bar{u}_i} \}$
that are in the image of $S_2$ can be attached
to $F_0$ along common boundary curves among $\{u_1, \dots u_n \}$,
forming a connected surface $S_3$ of genus $n$.
Finally, a connected surface $S_4$ with boundary $K$ is obtained
by identifying pairs of curves in $\partial S_3$ that are mapped
to the boundary components $c_{1}, \dots c_{m}$.
Each $c_i$ then has either two or four curves in $\partial S_3$ mapped to it, so there
are either one or two identifications made along each $c_i$.  The choice of which
pairs to identify, in case there are four curves mapped to $c_i$, is not important.
The Euler characteristic of $S_4$ is the same
as that of $S_3$ and $S_2$,
and therefore greater or equal to that of $S_1$.
So the genus of $S_4$ is at most that of $S_1$, genus$(S_4) \le n+m$.
Each identification of a pair of surface boundaries along
$c_{1}, \dots c_{m}$ contributes one to the genus of $S_4$.
There is at least one such identification
along each of the $m$ curves  $\{c_{j}\}$, though there
may be two if each of the three surfaces meeting
$F_0$ along $c_{j}$ has multiplicity one in $S_2$.
So identification of curves along $c_{1}, \dots c_{m}$
adds at least $m$ to the genus of $S_4$,
and it follows that genus($S_4) \ge n+m$. Since we have seen that
genus$(S_4) \le n+m$, equality must hold. Equality holds
when exactly one of the three surface pieces
meeting $F_0$ along $c_{j}$ has odd multiplicity for
each $1 \le j \le m$.
We then assign the value ``TRUE'' to a literal
$u_i$ if $F_{u_i}$ is used in $S_4$, and the value ``FALSE'' to $u_i$
if $F_{\bar{u}_i}$ is used in $S_4$.
This gives a truth assignment to $U$ in which each
clause in $C$ has exactly one true literal.
\qed

To show that {\sc 3-MANIFOLD KNOT GENUS} is {\bf NP}-hard
we reduce in polynomial time
an instance of  the  {\bf NP}-hard problem
{\sc ONE-IN-THREE SAT} to an instance
of  {\sc 3-MANIFOLD KNOT GENUS}.

\noindent
{\bf Proof of Theorem~\ref{thmnph}.}
Given an instance of  {\sc ONE-IN-THREE SAT} form $B$ as in
Lemma~\ref{branch}. Then there is a truth assignment for $U$ such 
that each clause in $C$
has exactly one true literal if and only if
$K$ is the boundary of a surface of genus $m+n$
mapped continuously into $B$.
Form a 3-manifold $N$ by thickening $B$ so that it is embedded inside a
triangulated 3-dimensional manifold with boundary. The thickening process
replaces each subsurface forming $B$ by a product of a surface with an
interval, and then glues these surfaces together along portions of
their boundaries, as indicated in Figure~\ref{thicken}.
The curve $K$ remains on the boundary
of $N$, and there is a projection map $p:N \to B$ which fixes $K$.
Form a closed manifold $M$ by doubling $N$ along its boundary, namely by
taking two copies of $N$ and identifying them along their boundaries
by the identity map.
Then $M$ admits an involution $\tau $ that fixes $K$, and has quotient $N$.

To find a triangulation of $N$, 
we first describe an explicit triangulation of the branching surface $B$.
An orientable surface with boundary has a triangulation with one
vertex on each boundary component and no vertices in the interior.
The number of triangles is $4g+5c-4$, where $g$ is the genus
and $c$ is the number of boundary components. We choose such triangulations
for each of the $2n+1$ subsurfaces in the branching surface, and we
match them together along boundary components to get a triangulation
of $B$.  So $B$ has a triangulation in which the number of
triangles is linearly bounded in $n+m$. We
thicken the surfaces $F_{u_i}$ and $F_{\bar{u_i}}$,
by taking their product with an interval.
We form a cell structure
of the thickened surface by dividing
the product into prisms, products of a triangle and an interval.
For  the thickening of $F_0$, we start by doing the same and then
go on to divide each interval into five subintervals.
The top, middle and bottom subintervals are identified with the
intervals from the three thickened surfaces meeting each thickened boundary
component of  $F_0$ corresponding to a curve $c_i$.  Only the top and bottom
intervals are identified with thickened surfaces from the other boundary
components.
We can therefore now glue these thickened surfaces together to get a
3-cell structure on $N$.  Each cell in this structure is
a prism.
We form a closed 3-manifold $M$ by doubling $N$ along its boundary,
gluing two copies of $N$ together
along their boundaries to obtain a 3-manifold with no boundary.

Finally, we stellar subdivide the
cell structure to get a triangulation. Each prism is
divided into 14 tetrahedra, by dividing each rectangular
face into four triangles by coning to a vertex in the
center of each such face, and then coning the 14 triangles
of the boundary of the prism to a vertex added to its center.
The number of resulting 
simplices in $M$ is linearly bounded by $n+m$. 

We now check that if $K$ bounds
a surface of genus  $g \le m+n$ with interior in
$M \backslash K$ then it bounds a surface of the same genus in $B$.
Suppose that $F$ is an embedded surface in $M$ with 
boundary $K$.  If
$F$ does not already lie in $N$, then perturb it slightly so
that the interior of $F$ meets $\partial N$
transversely in a finite number of simple closed curves and arcs.
Using the involution $\tau $, reflect the portion of $F$ not in
$N$ into $N$, forming an immersed surface $F'$
lying in $N$ and with the same boundary as $F$.
The interior of $F'$ remains disjoint
from $K$ since $K$ is fixed by the involution.
The projection $p(F')$ is a 
surface of genus $g$ mapped into $B$ with boundary $K$.
So if $K$ bounds an embedded
surface of genus at most $m+n$ in $M$ then
it bounds a surface of genus at most $m+n$ mapped into $B$.
The converse was shown in the proof of Lemma~\ref{branch}.
It then follows from Lemma~\ref{branch} that $K$ bounds an embedded
surface of genus at most $m+n$ in $M$ if and only
if $C$ is satisfiable with
each clause containing exactly one true literal.

The construction of $B, N$ and $M$ described above
each requires only a linear number of steps in the size of the instance of
{\sc ONE-IN-THREE SAT} with which we started, so that the reduction
requires polynomial time.
\qed

\section{Orbits of interval isometries}
\label{pseudogroup}

In this section we develop a combinatorial procedure
that will allow us to count the components of a normal surface.
The procedure computes the
number of orbits of a collection of $k$ isometries
between subintervals of an interval
$[1,N] \subset \ZZ$ in time polynomial in
$k \log N$. By lining up the intersections of a normal surface with the
edges of a triangulation, we obtain such subintervals.
Arcs of the normal surface on
the faces of the triangulation give rise to correspondences
of these intersection
points which are subinterval isometries.
We can then apply the algorithm developed here
to count the number of components of a  normal surface.

Assume that we have a set of integers
$\{1, 2, \dots , N \}$ and
a collection of bijections,
$g_i : [a_i, b_i] \to [c_i, d_i]$, $1 \le i \le k$,
either increasing or decreasing,
that are called {\em pairings}.
If a pairing identifies two intervals
$[a, b] $ and $[c, d ]$ by sending $a$ to $c$ and
$b$ to $d$, we call it an {\em orientation preserving} pairing, and if it
sends $a$ to $d$ and
$b$ to $c$, we call it {\em orientation reversing}.
If $a < c$ we refer to $[a,b]$ as the {\em domain} and
$[c,d]$ as the {\em range} of the pairing.
We work only with integers,
and use the term ``connected interval'' to refer to the integers in a
connected real interval.  The {\em width} of an interval $[a,b]$ with integer endpoints is
$b-a+1$, the number of integers it contains.
The {\em width} of a pairing is the width of its domain or range,
$w = b - a + 1 = d - c + 1 $.
If the pairing preserves orientation,
its {\em translation distance $t$} measures how far it
moves points, so $t = c - a = d - b$.
We can compose two pairings if the range of the first lies in the
domain or range of the second. The collection of pairings generates a pseudogroup,
under the operations of composition
where defined, inverses, and restriction to subintervals.

The interval $[1, \dots, N]$ is divided into equivalence classes by the
action of the pairings, which are called {\em orbits}.
We are interested in the orbit structure
of the collection of pairings, since with appropriate interpretation an orbit
corresponds to a connected component of a normal surface.
We introduce several simplification processes on the set of pairings in order
to analyze the structure of the set of orbits.

We introduce some terminology to describe the behavior of pairings.
An interval is called {\em static} if it is in neither the domain nor the
range of any pairing, so that its
points are identified to no other points by pairings.
Given a collection of pairings acting on the integers $[1, \dots, N]$, a pairing is said to 
be {\em maximal} if its range contains both $N$ and
the range of any other pairing containing $N$.
More precisely, define a linear order on pairings
using the lexicographical order $(d_i, -c_i, -a_i, -\mbox{orientation})$,
so that the maximal pairing has the highest upper endpoint, and
among those with that endpoint the widest range, and among those
with that range the biggest translation distance (if orientable).
Finally, we say that a pairing  $g : [a, b] \to [c, d]$ is {\em periodic} with period $t$ if
it is orientation preserving with translation distance $t$
and  $a < c = a+t \le b+1 $, so there is no gap between the domain and
range.  The combined interval $[a,  d]$ is then called a {\em periodic interval} of period $t$.

The following lemma describes the orbits of a periodic pairing.

\begin{lemma}
A periodic pairing has $t = c - a$ orbits on $[a,  d]$.
\end{lemma}
\bproof
Each point in $[a, b]$ greater or equal
to $c$ lies in the range of $g$ and
can be mapped to a smaller point in $[a, b]$ by a
power of $g^{-1}$.  So each orbit on $[a,  d]$ has a representative
in $[a, c - 1 ]$.  Since the congruence class modulo $t$ of
a point is preserved
by $g$, each of the $t$ integers in $[a, c - 1 ]$ lies in a distinct
orbit. These points uniquely represent the $t$ orbits.
\qed

We now show how to merge two pairings
with sufficient overlap into a single pairing with the same
orbits.

\begin{lemma}\label{merge1}
Let $R_1$ be a periodic interval with pairing $g_1$ of period $t_1$.
Suppose that there exists an orientation preserving
pairing $g_2$ with translation distance $t_2$
and an interval $J_1 \subset R_1$ such that
$J_1$ has width $t_1$ and $g_2(J_1) \subset R_1$. 
Then the orbits of $g_1 \cup g_2$ on $R_1$ are
the same as those of a single periodic action on $R_1$ of period
$ \mathrm{GCD}(t_1, t_2)$.
\end{lemma}
\bproof
Let $J_1$ be an interval in $R_1$ of width
$t_1$ that is in the domain of $g_2$ and which is paired
by $g_2$ to an interval in $R_1$. Each
point in $R_1$ has a unique orbit representative
in $J_1$ under the action of $g_1$ on $R_1$.
The interval $g_2(J_1)$ lies
in $R_1$ by assumption.  For $x \in J_1$ let
$f(x)$ be the unique point in $J_1$ obtained by
carrying $g_2(x)$ back to a point in $J_1$ by a power of $g_1$ or  $g_1^{-1}$.
The effect of $f(x)$ on $J_1$ is a shift of $t_2 \mbox{ (mod }t_1)$.
The orbits of $f$ on $J_1$ divide the points of $J_1$ into
congruence classes modulo $\mathrm{GCD}(t_1, t_2)$.
Neither $g_1$ nor
$g_2$ change the congruence class of a point mod $ (\mathrm{GCD}(t_1, t_2))$,
so a subinterval of width $ \mathrm{GCD}(t_1, t_2)$ in
$R_1$ contains exactly one representative of each orbit
of $R_1$ under the action of $g_1 \cup g_2$. The same orbits
arise from a periodic action of period
$ \mathrm{GCD}(t_1, t_2)$.
\qed

The following is a special case of Lemma~\ref{merge1}
that applies when both $g_1$ and $g_2$ are periodic pairings.

\begin{lemma}\label{merge2}
Let $R_1$ and $R_2$ be overlapping periodic intervals,
associated with pairings $g_1, g_2$ having periods $t_1$ and $t_2$.
Suppose that width($R_1 \cap R_2) \ge t_1 + t_2$.
Then the orbits of  $g_1 \cup g_2$ on $R_1 \cup R_2$ are
the same as those of a single periodic
pairing on  $R_1 \cup R_2$ of period $ \mathrm{GCD}(t_1, t_2)$.
\end{lemma}
\bproof
The leftmost interval $J_1$ of width $t_1$ in $R_1 \cap R_2$
is translated by $g_2$ a distance of $t_2$ to the right,
to an interval $J_2$ which lies in $R_1 \cap R_2$.
Lemma~\ref{merge1} then states that the action of  $g_1 \cup g_2$ on $R_1$ is
the same as those of a single periodic
pairing on $R_1$ of period $ \mathrm{GCD}(t_1, t_2)$.
By symmetry the same result holds for the action of 
 $g_1 \cup g_2$ on $R_2$. Thus the  orbits of  $g_1 \cup g_2$ on $R_1 \cup R_2$ are
the same as those of a single periodic
pairing on  $R_1 \cup R_2$ of period $ \mathrm{GCD}(t_1, t_2)$.
\qed

We state a consequence in a form which will
be convenient for our applications.

\begin{lemma}\label{merge3}
Let $g_1, g_2$ be periodic pairings with periods $t_1, t_2$ and let $J_1, J, J_2$
be intervals with
$ g_1(J_1) = J $ and $g_2(J) = J_2 $.
Suppose that $J$ is contained in the union $ J_1 \cup J_2$.
Then the hypothesis of Lemma~\ref{merge2} is satisfied and
the orbits of $g_1 \cup g_2 $ are
the same as those of a single periodic pairing of period
$\mathrm{GCD}(t_1, t_2)$, acting on the union of the periodic intervals
of $g_1, g_2$.
\end{lemma}
\bproof
Let $ J_1 = [a_1,b_1], J= [a,b]$ and $J_2 = [a_2,b_2]$.
We have $a_2 \ge a_1$ and $J \subset J_1 \cup J_2 $.
Since $J_1 \cup J_2$ is a connected interval, we have that $a_2 \le b_1 + 1 $.
Note that $t_1 =  b - b_1 = a - a_1$ and $t_2 = b_2 - b = a_2 - a$.
Then $t_1+t_2 = (b - b_1) + (a_2 - a) = (b-a+1) + (a_2 - b_1 - 1)$.
The width of $J$ is $(b-a+1)$ and we have seen above
that $(a_2 - b_1 - 1) \le 0$.  So the width of $J$ is at least
$t_1+t_2$ and Lemma~\ref{merge2} implies the conclusion of the lemma.
\qed

The orbit counting algorithm applies
a series of modifications to a collection of pairings.
We now describe these modifications.
\begin{itemize}
\item Periodic merger

The  {\em periodic merger} operation replaces $g_1$ and $g_2$ by a single periodic
action on  $R_1 \cup R_2$ of period $ \mathrm{GCD}(t_1, t_2)$, as in
Lemma~\ref{merge2}.

\item Contraction

The operation of {\em contraction} is performed
on a static interval $[r,s]$.
We eliminate this interval, replace $[1,N]$ by
$[1, N-(s-r+1)]$, and alter each $g_j$ by replacing any point $x$ in a
domain or range which lies entirely to the right of $s$ by $x-(s-r+1)$.
(This operation will lead to a decrease of $s-r+1$ in the
number of orbits, since the eliminated points are
each unique representatives of an orbit.)

\item Trimming

The {\em trimming} operation simplifies an orientation reversing pairing whose
domain and range overlap.
Suppose that $g : [a, b] \to [c, d]$ is a pairing
with $g(a) = d, g(b) = c$ and $b \ge c$.
Define a new pairing $g': [a, (a+d)/2 ) \to ((a+d)/2, d]$
by restricting the domain and range of $g$,
and say that $g'$ is obtained from $g$ by {\em trimming}. The domain and
range of a trimmed pairing are disjoint.

\item Truncation

If an interval lies in the domain and range of exactly one
pairing, then the interval can be ``peeled off'' without changing the
orbit structure, in a way we now describe.  The algorithm applies this operation
to strip off points from the right of the interval $[1,N]$.

When there is a pairing $g : [a, b] \to [c, N]$ and a
value $N'$ with $c \le N'+1 \le N$, such that all points
in the interval $[N'+1,N]$ are in the
range of no pairing other than $g$, then
we can perform an operation called {\em truncation} of $g$.
Truncation shortens the interval $[1, N]$ to the interval $[1,N']$,
and similarly shortens the domain and range of $g$.
If $g$ is orientation preserving,
pairings other than $g$ are unchanged, while
$g$ is eliminated entirely if $c=N'+1$, or
replaced by a shortened pairing $g': [a, b - (N-N')] \to [c, N']$
if $c \le N'$.
We can perform this operation even if the interval $[N'+1, N]$
intersects both the range and the domain of $g$.

If $g$ is orientation reversing,
truncation is applied only when $g$ has disjoint domain
and range (i.e. after trimming), and $[N'+1,N]$ is contained in the range.
Suppose $g: [a, b] \to [c, N]$ is a pairing
with $g(a) = N, g(b) = c$ and $b < c$, that
$[N'+1,N]$ is disjoint from the domains and ranges
of all pairings other then $g$ and that $c \le N'+1 $.
If $N'+1=c$ then we  eliminate $g$.
Otherwise replace $g$ by
a shortened orientation reversing pairing, formed by restricting its
domain to $ [a+N-N', b]$ and its range to $[c, N']$.

\item Transmission

Transmission is an operation in which two pairings are composed.
A pairing $g_1$ used to shift down the domain
and range of a second pairing $g_2$ as much as possible.
This operation will allow us to shift pairings leftwards from the right end of the
interval $[1,N']$ and subsequently apply truncation.

If $g_1$ is orientation reversing and has overlapping
domain and range, then as a first step in transmission we trim $g_1$.
Now consider a pairing $g_1$, either orientation preserving or orientation reversing,
and a second pairing $g_2$
whose range is contained in the range of $g_1$.
If the domain of $g_2$ is not contained in the range of $g_1$,
form the composite map
$g_2' = g_1^{-r} \circ g_2$, where $r=1$ if $g_1$ is orientation reversing
and otherwise $r \ge 1$ is the largest integer
such that $g_1^{-r +1} ( [c_2, d_2])$
is contained in the range of $g_1$.
The domain of $g_2'$ is the same as that of $g_2$ in this case.
If the domain of $g_2$ is also contained in the range
of $g_1$ then form the composite map
$g_2' = g_1^{-r} \circ g_2 \circ g_1^s : g_1^{-s}([a_2, b_2]) \to g_1^{-r}([c_2, d_2])$,
where $r$ is as above,
$s=1$ if $g_1$ is orientation reversing, and
otherwise $s \ge 1$ is the largest integer such that
$g_1^{-s+1} ( [a_2, b_2] )$
is contained in the range of $g_1$.
The domain of $g_2'$ is then that of $g_2$ shifted left by $g_1^{-s}$.
The process of replacing $g_2$ by
$g_1^{-r} \circ g_2 \circ g_1^s $ is called a
{\em transmission of $g_2$ by $g_1$}.

\end{itemize}

We now construct a sequence of pseudogroups of pairings, terminating with the trivial
pseudogroup acting on the empty set.
At each stage, each orbit is associated to a unique orbit in the previous set of pairings,
though the number of orbits may decrease. A counter is kept at each stage
that records the total decrease in the number of orbits.  Since the final pseudogroup 
has no orbits, the final value of this counter
gives the initial number of orbits.

\begin{lemma} \label{orbit.same}
A contraction decreases the number of orbits by the width
of the contracted interval.
Altering a collection of pairings by any of the
operations of periodic merger, trimming, truncation and transmission
preserves the number of orbits.
\end{lemma}
\bproof
We refer to the collection of pairings
$g_1 ,g_2, \dots g_k$ as $G$ and to the
new collection of pairings produced by one of the operations as $G'$.

The effect of a contraction is to 
shorten the interval $[1,N]$ by removing points
which are fixed by the entire collection of pairings. The number of orbits removed
equals to the width of the contracted interval.

A periodic merger joins two periodic intervals $R_1$ and $ R_2$ and their pairings
into one.
Lemma~\ref{merge2} shows that the orbits
in $R_1 \cup R_2$ of $g_1 \cup g_2$ are
the same as those produced by a single periodic pairing $g'$.
Suppose that $x$ and $y$ are points in  $[1,N]$ in the same orbit of
the action of $ G = \{g_1,g_2, \dots g_k\}$.
Then $h \cdot x = y$, where
$h$ is some finite word in the elements of $G$ and their inverses.
Wherever a $g_1$ or $g_2$ occurs in this word we
can replace it by a power of $g'$, since $g'$ translates
by an amount that divides the translation distance of $g_1$
and $g_2$. So $x$ and $y$ are in the same orbit of
$G' = \{ g', g'_3, \dots g'_k \}$. Conversely
suppose that that $x$ and $y$ are in the same orbit under
$G'$. Then  $h' \cdot x = y$, where
$h'$ is some finite word in $g', g'_3, g'_4 \dots g'_k$.
Lemma~\ref{merge2} implies that wherever a $g'$ occurs in
$h'$ we can replace it by some word in $g_1$ and $g_2$, since
the orbits of $g'$ and $g_1 \cup g_2$ coincide on the periodic interval
of $g'$. So periodic mergers preserve orbits.

Suppose next that $g_i: [a_i, b_i] \to [c_i, d_i]$ is an orientation
reversing pairing with $g_i(a_i) = d_i, g(b_i) = c_i \le b_i$,
and that
$g_i' :  [a_i,  (a_i+d_i)/2 )  \to ((a_i+d_i)/2), d_i]$
is obtained by trimming $g_i$.
If $x$ and $y$ are in the same orbit
of $ G = \{g_1,g_2, \dots g_k\}$,
then there is a sequence of points $x=x_1,x_2, \dots , x_r = y$ where
each $x_j$ is the image of $x_{j-1}$ under a pairing in $G$.
We can replace an occurrence of $g_i$ or $g_i^{-1}$ by
$g'$ if the point $x_{j-1}$ is smaller than
$(a_i+d_i)/2$ and by ${g'}^{-1}$ otherwise. So $x$ and $y$ remain in
the same orbit under the action of $G'$. Conversely
if $x$ and $y$ are in the same orbit under
$G'$ then replacing each occurrence of $g'$ by $g$ and 
of ${g'}^{-1}$ by $g^{-1}$ gives a word
in $G$ taking $x$ to $y$. We conclude that trimming preserves orbits.

Next consider the effect of a truncation.
Suppose that all points in $[N'+1,N]$ are in the
range of exactly one pairing
$g : [a, b] \to [c, N]$, and we truncate, shortening
to an interval $[1, N']$ and either eliminating $g$ or shortening it to
$g': [a, b - (N-N')] \to [c, N']$ in the orientation preserving case,
or to $g': [a + N-N', b] \to [c, N']$ in the orientation 
reversing case. Suppose that $y = h \cdot x$ where
$h$ is a reduced product of pairings (a product that does not contain
a subproduct of the form $g_j \cdot g_j^{-1}$).
We can assume that $y>x$, as otherwise we
can replace $h$ by its inverse.
Let $h'$ be obtained from $h$ by replacing all occurrences of $g$ with $g'$.
The successive images of $x$ under the subwords of $h$
can never enter and leave the interval $[N'+1,N]$, since
they can only enter it under the action of some positive power
of $g$, and only leave it under a negative power of $g$.
These don't occur in succession in a reduced product.
The image of a point $z$ under $g_j$ is unchanged
unless $g_j=g$, and $g(z) \ge N'+1$.
So the image of $x$ under $h$ is the same as its image
under $h'$ unless $y > N'$. In that case $h$ is of the
form $h = g_i^r h_1$ where $h_1 \cdot x \le N'$, $r \ge 1$,
and $ y = g_i^r h_1 \cdot x > N'$.
So two points in $[1, N']$ are in the same orbit under
$G'$ if and only if they are in the same orbit under $G$.
But each point in $[N'+1,N]$ is in the same $G$ orbit as a point in $[1,N']$.
It follows that the number of orbits is unchanged by truncation.

Finally suppose that a transmission of $g_j$ by $g_i$
replaces the pairing $g_j$ in $G$ by $g_j' = g_i^{-r} \circ g_j \circ g_i^s$.
If $y = h \cdot x$ where $h$ is a word in $G$, then
$y = h' \cdot x$ where $h'$ is obtained by replacing every occurrence of
$g_j$ with  $g_i^k \circ g_j \circ g_i^{-s}$. Similarly if $y = h'
\cdot x$ where $h'$ is a word in $G'$, then
$y = h \cdot x$ where $h$ is obtained by replacing every occurrence of
$g_j'$ with  $g_i^{-k} \circ g_j \circ g_i^s$.
So transmissions preserve the collection of orbits and don't change
their number.
\qed

We now describe an algorithm which uses these operations to count the number
of orbits of a pseudogroup of pairings acting on $[1,N]$.
We initially set a counter for the number of orbits to zero.
A pairing $g_i$ is said to be {\em maximal} if
$d_i = N$ and $[c_i,N]$ contains the range of any other pairing with
an endpoint at $N$.
Let $S = \{ 1, 2, \dots N \}$ and let $g_i : [a_i, b_i] \to [c_i,
d_i],\  1 \le i \le k$
be a collection of pairings between subintervals of $S$.
The algorithm will contract $S$ and reduce the number of pairings.
It keeps a running count of the number of orbits detected in an integer
referred to as the orbit counter.
Denote by $N'$ the current size of the interval as we proceed.
The algorithm repeats
the following steps, reducing $N'$ and $k$, until there are no
points remaining.\\

\noindent
{\bf Orbit counting algorithm:}
\begin{enumerate}
\item 
Delete any pairings that are restrictions
of the identity.
\item 
Make any possible contractions and, if any exist, increment the orbit counter by the 
sum of the number of points deleted by the contractions.
If the number of pairings remaining is zero, output the number of orbits and stop.
\item 
Trim all orientation reversing pairings whose domain and range
overlap.
\item 
Search for pairs of periodic pairings 
$g_i$ and $g_j$ whose domains and ranges
satisfy the condition of  Lemma~\ref{merge2}.
If any such pair exists
then perform a merger as in Lemma~\ref{merge2},
replacing $g_i$ and $g_j$ by a single periodic pairing,
with translation distance $\mathrm{GCD}(t_i, t_j)$.
The new pairing acts on the union of the domains and ranges of $g_i$ and $g_j$.
Repeat until no mergers can be performed.
\item 
Find a maximal $g_i$. For each $g_j \ne g_i$ whose range
is contained in $[c_i, N']$, transmit $g_j$ by $g_i$.
\item 
Find the smallest value of $c$ such that the interval $[c, N']$
intersects the range of exactly one pairing.
Truncate the pairing whose range contains the interval $[c, N']$.
\end{enumerate}

\begin{theorem} \label{orbitcount}
The orbit counting algorithm gives the number of orbits of
the action of the pairings $\{g_i\}_{i=1}^k$ on $[1,\dots ,N]$ in time
bounded by a polynomial in $ k \log N $.
\end{theorem}
\bproof
We first check that the orbit counting algorithm correctly counts the number
of orbits.
In Step (1), deleting a pairing which is the identity on its domain does
not change the number of orbits.
In Step (2), contracting a static interval 
removes a number of points that are unique orbit representatives,
and the count of these is added to the running total kept in the orbit counter.
If all the points have been removed, then there are no more orbits
to count and the orbit count is complete.
Lemma~\ref{orbit.same} shows that the operations of
transmission, trimming, truncation
and merger occurring in Steps (3),(5) and (6)
do not change the number of orbits of the collection of
pairings.

Mergers carried out in Step (4)
occur when the conditions of  Lemma~\ref{merge2}
are satisfied, and  Lemma~\ref{merge2}
shows that they
preserve the orbit structure of the collection of pairings.
A merger reduces the
number of pairings by one, replacing two pairings $g_i$ and $g_j$
by a single periodic pairing acting on the union of their periodic intervals.

In each cycle through these steps
the interval width decreases by at least one, either in
Step (2) or in Step(6).
It follows that the algorithm terminates,
yielding a count of the number of orbits
after a number of steps bounded by
the width $N$ of the interval.
We will obtain a much better bound.
To do so, we define a complexity which decreases as we
iterate the above steps.
Recall that $w_i = b_i - a_i + 1,\  1 \le i \le k$ and
$k$ is the number of pairings. The complexity $X$ is defined to be
$$
X = 4^k \prod_{i=1}^k w_i .
$$
The process of executing Steps (1)-(7) in turn is called a {\em cycle}.
We will show that when we run through $2k$ cycles,
$X$ is reduced by a factor of at least two.
See the remark at the end of this section
for a geometric interpretation of this complexity.

Call an interval $[x,y]$ a {\em $Z$-close} interval
if it is the domain or range of a pairing of
subintervals of $[1,Z]$ and if $ y-x+1 \ge (Z-x+1)/2$, or
equivalently, if $y \ge Z/2 +x/2 -1/2$.
Being $Z$-close corresponds to being relatively close
to $Z$. It means that there is no room for
another interval of the same size to the right of $y$.
The value of $Z$ is initially set to $N$,
and as the algorithm proceeds it 
is reset to the current interval width $N'$ each time
the number of pairings $k$ is decreased. The number of pairings
decreases when a merger occurs, or when a pairing is
truncated to zero width, or when a pairing is transmitted
to become trivial (a restriction of the identity) and then
eliminated.

\begin{claim} \label{overlap}
The union of two $Z$-close intervals of equal
width is a connected interval.
\end{claim}
\bproof
Suppose that $[a,b]$ and $[c,d]$ are equal width $Z$-close intervals
with $c \ge a$.
Then $b-a = d-c$, $b-a+1 \ge (Z-a+1)/2$ and $d-c+1 \ge (Z-c+1)/2$,
so $b \ge Z/2 +a/2 -1/2 \ge d/2 +a/2 -1/2
= b/2 + c/2 -1/2 $ implying that $b \ge c-1$ and
$[a,b] \cup [c,d] $ is connected.
\qed

\begin{claim} \label{transmit}
Suppose that the domain $[a,b] \subset [1,Z]$ of a pairing $g$ is not $Z$-close.
Then the image $[x',y']$ of an interval $[x,y]$ under  ${g}^{-1}$ is 
not $Z$-close.
\end{claim}
\bproof
Suppose that  $[x',y']$ is the image of $[x,y]$ under
${g}^{-1}$, where $g: [a,b] \to [c, d]$.
By assumption we have that $b < Z/2 + a/2 -1/2$.
The interval $[x',y']$ is contained in $[a,b]$ so $a \le x' \le 
y' \le b$.
Then $y' < Z/2 + a/2 -1/2 \le  Z/2 + x'/2 -1/2$ and
$[x',y']$ is not $Z$-close.
\qed

\begin{claim} \label{5pairs}
After a series of five cycles
either the number of $Z$-close intervals decreases or the number of
pairings decreases.
\end{claim}
\bproof
In each cycle at least one maximal pairing is truncated.
Suppose that during a series of five cycles, five successive maximal pairings
$g_{1}, g_{2}, g_{3},g_4,g_5$ occur in turn in the orbit counting algorithm.
The initial maximal pairing $g_1$ is truncated in the first
cycle. Eventually $g_{1}$ stops being maximal, and it is then
transmitted by the new maximal pairing $g_{2}$ in the next cycle.
The cycle $g_2$ is then itself truncated until no longer maximal and it
in turn is then transmitted by $g_{3}$ and so on. When
$g_{i+1}$ transmits $g_{i}$, the pairing of 
$[a_{i}, b_{i}]$ with $[c_{i}, d_{i}]$ is replaced
by a pairing of $[a_{i}', b_{i}']$ with $[c_{i}', d_{i}']$ and 
one of the following three cases occurs.

\begin{enumerate}
\item 
The range $[c_i, d_i]$ of $g_i$
is transmitted to a non $Z$-close interval by $g_{i+1}$.
\item 
The range $[c_i, d_i]$ is transmitted to a $Z$-close interval
$[c_i', d_i']$ by $g_{i+1}$, and the domain $[a_i', b_i']$ of the transmitted pairing is $Z$-close.
\item 
The range $[c_i, d_i]$ is transmitted to a $Z$-close interval
$[c_i', d_i']$ by $g_{i+1}$, and the domain $[a_i', b_i']$ of the transmitted pairing is not $Z$-close.
\end{enumerate}

We first consider the case where each of $g_1, g_2, and g_3$ is orientation preserving.
Setting $Z$ to the initial interval width $N'$ and noting that $g_{1}$ is initially maximal, we see that
the range of $g_1$ is initially $Z$-close, with $Z = N' = d_1$. Truncation reduces the range to
of $g_1$ to $[c_{1}, d_{2}]$, at which point the pairing $g_{2}$
becomes maximal. We can assume the range $[c_{1}, d_{2}]$ of $g_1$ is still $Z$-close, or we
are done. In the next cycle the interval $[c_{1}, d_{2}]$ is transmitted by $g_{2}$
to an interval $[c_{1}', d_{1}']$.
If Case (1) applies, this new interval is not $Z$-close and the number of
$Z$-close intervals has decreased.
If Case (2), the three intervals $[c_1, d_2], [c_1', d_1']$ and
$[a_1', b_1']$ are all $Z$-close and of equal width, and
therefore satisfy the hypothesis of Claim~\ref{overlap}.
Lemma~\ref{merge3} then implies that
two pairings can be merged, and the number of pairings
decreases during the execution of Step (4) in the next cycle.
If Case (3) occurs, the maximal pairing $g_{2}$ is truncated in the second cycle.
We can assume its domain remains $Z$-close, or we are done.
When the next pairing $g_{3}$ becomes maximal and transmits $g_2$ in the third cycle,
we must fall into one of the first two cases. So
the number of $Z$-close
intervals is decreased in the third cycle
or the number of pairings decreases in a merger during a subsequent application of Step (4)
in the fourth cycle.

Now we consider the possibility that one of  $g_1, g_2, and g_3$ is orientation reversing.
The domain and range of a trimmed orientation reversing
pairing are disjoint, and hence any interval transmitted by an
orientation reversing pairing is not  $Z$-close.
If any of $g_2, g_3 or g_4$ are  orientation reversing 
then it transmits the previously maximal pairing's range to
an interval that is not  $Z$-close and the number of  $Z$-close pairings decreases.
If non is orientation reversing then
there is a sequence of three successive orientation preserving
pairings and the previous argument applies.
\qed

\begin{claim} \label{div2}
Suppose that $Z$ is set to the
current interval size $N'$ and that
after a series of truncations in which the maximal pairings
are successively $g_1, g_2, \dots ,g_r$,
no $Z$-close intervals remain. Then the complexity $X$ is reduced
by a factor of at least two.
\end{claim}
\bproof
Truncation of a maximal pairing
$g_i$ results in its range $[c_i,d_i]$, of width $w_i$, being reduced
to a shorter interval $[c_i, d_i']$ of width $w_i'$.
This reduces the complexity $X$ by a factor of
$$
\frac{w_i'}{w_i} = \frac{d_i' - c_i +1}{d_i - c_i +1}.
$$
The maximal pairing changes after a truncation if the range of
$g_i'$ has truncated sufficiently so that it
is contained in the range of $g_{i+1}$, so that
$c_{i+1} \le c_i'$ and $ d_{i+1} = d_i'$.
Define $d_i''$ by 
$$
d_i'' = c_i -1 +  (w_{i+1}) \frac{Z-c_i + 1}{Z-c_{i+1} + 1}.
$$
and define
$$
w_i'' = d_i'' - c_i +1.
$$
Differentiation shows that the function
$$
f(x) =  \frac{x - c_i +1}{x - c_{i+1} +1}
$$
is increasing with $x$ for $x \ge  c_i$, since $c_i \ge c_{i+1} $. So
\begin{eqnarray*}
w_i' & = & d_i' - c_i + 1 = d_{i+1} - c_i + 1 =
(w_{i+1}) \frac{(d_{i+1} - c_i + 1)}{w_{i+1}}
\\
& = &
(w_{i+1}) \frac{(d_{i+1} - c_i + 1)}{(d_{i+1} - c_{i+1} + 1)}
\le (w_{i+1})  \frac{Z - c_i +1}{Z - c_{i+1} +1} \\
 & = & d_i'' - c_i + 1 = w_i'' 
\end{eqnarray*}

After a series of truncations of the pairings
$g_1, g_2, \dots  ,g_r$, the complexity $X$
is multiplied by a factor of
\begin{eqnarray*}
&    &  \left(\frac{w_1'}{w_1}\right)
     \left(\frac{w_2'} {w_2}\right) \dots
     \left(\frac{w_r'} {w_r}\right)  \\
& \le & \left(\frac{w_1''}{w_1}\right)
     \left(\frac{w_2''} {w_2}\right) \dots
     \left(\frac{w_{r-1}''} {w_{r-1}}\right)
     \left(\frac{w_r'} {w_r}\right)  \\
& = & \left( \frac{d_1'' - c_1 +1}{w_1}\right)
     \left( \frac{d_2'' - c_2 +1}{w_2}\right) \dots
     \left(\frac{d_{r-1}''  - c_{r-1} +1} {w_{r-1}}\right)
     \left( \frac{w_r'} {w_r} \right) \\
& = &  \left( \frac{w_{2}}{w_1}\right)
      \left( \frac{Z-c_1 +1}{Z-c_{2} +1} \right)
    \dots
      \left( \frac{w_r}{w_{r-1}}\right)
      \left( \frac{Z- c_{r-1}+1}{Z-c_r +1}   \right)
        \left(\frac{w_r'} {w_r}\right) \\
& =  &  \left( \frac{Z-c_1 +1}{Z-c_r +1}
\right) \left( \frac{w_r'}{w_1} \right).
\end{eqnarray*}
Now  $Z = d_1$, since $g_1$ was the first maximal pairing
truncated, so $Z-c_1 +1 =   d_1  -c_1 +1 = w_1$.
By assumption, there are no
$Z$-close intervals remaining after $g_r$ is truncated,
so that $[c_r, d_r']$ is not $Z$-close. We then have by the
definition of $Z$-close that
$$
\left( \frac{Z-c_1 +1}{Z-c_r +1}  \right)
\left( \frac{w_r'}{w_1} \right) =
\frac{w_r'}{Z - c_r + 1} = \frac{d_r' - c_r + 1}{Z -
c_r + 1} < 1/2 .
$$
It follows that $X$ is multiplied by a factor smaller than 1/2.
\qed

\noindent
{\bf Proof of Theorem~\ref{orbitcount}:}
The operations involved in one cycle consist
of comparisons, additions, subtractions
and computing greatest common divisors of the $k$ pairings  of $G$.
The number of these operations occurring in each cycle
is linear in $k^2 \log N$.
Pairings are described by pairs of intervals,
whose boundary points are integers of size at most $N$. So the 
running time of each cycle is polynomial in $k \log N  $.

If the number of pairings $k$ decreases during a cycle,
then a pairing has been eliminated because its width has
truncated to zero, because it has been transmitted and become
the identity on its domain, or because two pairings have merged.
When two pairings $g_1, g_2$ merge, the value of $k$ decreases by one.
The product $(w_1) (w_2)$ which occurs in $X$
is replaced by the width of the
new pairing, which is at most $(w_1 + w_2) \le 2w_1w_2$.
Therefore $X = 4^k \prod_{i=1}^k w_i$ decreases by a factor of
at least two after a merger.

Since each pairing has width at most $N$,
the initial complexity is bounded above by $N^k 4^k$. 
By Claim~\ref{5pairs}, each time we run through five cycles
either the number of
pairings decreases or the number of $Z$-close intervals decreases.
There are at most $k$ $Z$-close intervals initially, so after
$5k$ cycles either there is a reduction in the number
of pairings or there are no $Z$-close intervals
remaining.  In the second case the complexity $X$ has decreased by a
factor of at least two, by Claim~\ref{div2}.

So $X$ decreases by a factor of at least two after $5k$ cycles,
and the complexity reduces to zero
after at most $k (2 + \log_2 N) $ successive series of $5k$ cycles, or
after at most $5k^2(2 + \log_2 N) $ cycles.
Each cycle runs in time polynomial in
$k \log N$, so the total running time
is also polynomial in $k \log N$.
\qed

We now apply Theorem~\ref{orbitcount} to count the number of components of
a normal curve or normal surface.  An obvious algorithm to count components
proceeds by marking vertices connected by
common edges until all vertices in a component are reached.
This procedure takes time linear
in the number of edges of the curve. This is equal to
the sum of the normal coordinates (called $W$ below), but we can achieve an
exponential improvement.
We first look at normal curves.

\begin{figure}[hbtp]
\centering
\includegraphics[width=.6\textwidth]{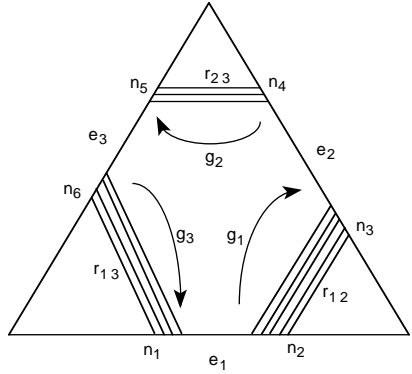}
\caption{ \label{pairings}
Normal arcs on a triangular face give three pairings of
ordered intervals:
$g_1 : [n_1+r_{13}, n_2] \to [n_3, n_3+r_{12}-1]$,
$g_2 : [n_3+r_{12}, n_4] \to [n_5, n_5+r_{23}-1]$, and
$g_3 : [n_5+ r_{23}, n_6] \to [n_1, n_1+r_{13}-1]$
The normal curve coordinates associated
to this face are $r_{12}, r_{13}$ and $ r_{23}$.
In a normal surface, $r_{12}, r_{13}$ and $ r_{23}$
are each a sum of two normal coordinates.}
\end{figure}

\begin{corollary} \label{normalcurvecount}
Let $F$ be a surface with a triangulation $T$
containing $t$ triangles and
let $\gamma$ be a normal curve in $F$
with normal coordinates summing to $W$.
There is a procedure for counting the number of components of $\gamma$
that runs in time polynomial in $t \log W$.
\end{corollary}
\bproof
The  1-skeleton of $T$ contains $e$ edges, where $e \le 3t$.
Fix once and for all an ordering of these edges.
A normal curve $\gamma$ intersects each edge of $T$ in 
a finite number of points.
Set $N$ to be the weight $W$ of $\gamma$, the
sum of the number of intersection points of $\gamma$ with all the edges of $T$.
Label the intersections of $\gamma$ and the first edge
of the 1-skeleton by the integers $1, 2, \dots, i_1$, the intersections
of $\gamma$ and the $j$th edge by $i_{j-1}+1, i_{j-1}+2, \dots, i_j$, and the
intersection of $\gamma$ and the $e^{th}$
edge of $T$ by $i_{e-1}+1, i_{e-1}+2, \dots, i_e$.
Then $i_e = N$. 
Each triangular face of $T$ has three sets of arcs pairing points of $[1,N]$,
with one set running between each pair of edges of the face.
To each set of arcs we associate a pairing between the
intervals at either end of the arcs, as in Figure~\ref{pairings}.
All the pairings are orientation reversing in this example.
In general some will be orientation preserving, as
the edge orientations on any triangle can be arbitrary.

The number of connected components of
$\gamma$ is the same 
as the number of orbits of  an action of a collection of pairings on
$\gamma \cap T^{(1)}$,
where two points are paired if they are connected by
an edge of $\gamma $ lying in a triangle.
This is precisely the number returned by the orbit counting algorithm.
The number of pairings is at most $3t$ and the interval size is $W$.
Applying Theorem~\ref{orbitcount}, 
we can determine the number of components of
the normal curve $\gamma$
in time polynomial in $t \log W$.
\qed

A similar argument applies to normal surfaces.

\begin{corollary} \label{normalcount}
Let $M$ be a 3-manifold with a triangulation $T$
containing $t$ tetrahedra and
let $F$ be a normal surface in $M$ whose normal
coordinates sum to $W$.
There is a procedure for counting the number of components of $F$
which runs in time polynomial in $t \log W$.
\end{corollary}
\bproof
The 2-skeleton of $T$ contains at most $4t$ faces and the
1-skeleton contains $e\le 6t$ edges.
Set $N$ to equal the weight $W$ of $F$, the total number of points
in which it intersects the 1-skeleton.
Order the edges of $ T^{(1)}$ in an arbitrary way and
label the intersections of $F$ and the
$j$th edge of  $ T^{(1)}$ by $i_{j-1}+1, i_{j-1}+2, \dots, i_j$.
Again we have $i_e = N$. 
To a pair of edges on a triangular face of the 2-skeleton
we associate a pairing between the
intervals at either end of the corresponding set of arcs.
There are at most three pairings for each
face, and the number of faces is at most $4t$,
so the number of pairings is bounded above by $12t$.
These pairings are determined by the normal coordinates of $F$.

For a normal surface, the number of connected components of
$F$ is the same as the number of components of
$F \cap T^{(2)}$, since every component of $F$ intersects the 1-skeleton of $T$,
and if two points on $F  \cap T^{(2)} \cap T^{(1)}$ can be joined by a path then
that path can be homotoped into the 2-skeleton of the triangulation.
So the the number of components equals
the number of orbits of $F \cap T^{(1)}$
under pairings that identify two points connected by
an edge of $F \cap T^{(2)}$ contained in a face of the 2-skeleton. 
This number is precisely what is computed
by the orbit counting algorithm.
Therefore we can determine the number of components of
the normal surface $F$ in time polynomial in $t \log W$.
\qed

\noindent
{\bf Remark:}
There is a motivating geometrical construction behind the combinatorics
of the orbit counting algorithm which we informally explain.  
We can associate to each pairing $g_i : [a_i, b_i] \to [c_i, d_i]$
a pair of ``transmission towers''
in the upper half-plane, with one tower being the vertical line segment
from $(a_i, 0)$ to $(a_i, w_i)$ and the other the vertical segment
from $(c_i, 0)$ to $(c_i, w_i)$.
These towers capture the information contained in the pairings.
Points in the domain and range of a pairing 
beam up to the tower leftward at a 45 degree angle, then beam
to the paired tower (either straight across, if the pairing
preserves orientation, or crossing if not) then down again.  Assign a ``cost''
to each transmission tower equal to the hyperbolic length in the upper
half space model from $y=1/e$ to the top of the tower.
This cost equals $\log (d_i - c_i+1 ) + 1$, and the sum of all these costs
is essentially the logarithm of the complexity $X$ used above.

The counting algorithm starts with
a Euclidean line emerging from the right endpoint of a $Z$-close interval
and going upward at a 45 degree angle. 
In hyperbolic terms, this is an equidistant curve from the geodesic
from the right endpoint to the point at infinity. 
Initially, we may assume that this sweep line hits the top of at least
one transmission tower. We make the highest of these the $Z$-close tower,
and use it to beam all the other towers to its paired interval.
If the two domains overlap, we use the highest power of the
transmission that can be applied. Eventually the triangle between the
current tower and the sweep line is vacant.
We sweep leftward with this equidistant curve.  As long as the sweep line
hits only one tower, we truncate it.

Consider a second equidistant curve $P$ through a point
$x$ on $\RR$, with slope -1/2.
The hyperbolic distance between the slope -1 and slope -1/2 lines is 
a constant equal to about 0.49.
Call the region between these two lines a {\em zone}.
If two paired towers 
have their tops in a single zone, their domains overlap.
If no towers are completely truncated away
between the time when the sweep line hits a point $x$
and when there are no towers intersecting the zone above $P$,
each equidistant curve between the two
given ones hits a transmission tower based at a point that will be removed by truncation.
So in this time either the number of towers has
decreased or the sum of all the tower costs has decreased by
at least 0.49.

If we merge transmission towers when possible,
a complexity based on the cost decreases
sufficiently fast to give a polynomial time algorithm.
The calculations in Theorem~\ref{orbitcount} implemented this geometric picture.

\section{{\sc 3-MANIFOLD KNOT GENUS} is {\bf NP}. }\label{sectiongnp}

In this section we establish that {\sc 3-MANIFOLD KNOT GENUS} is {\bf NP}.

\noindent
{\bf Proof of Theorem~\ref{thmgnp}:}
We begin with a simplicial complex consisting of
$t$ tetrahedra whose faces are identified in pairs, and a collection $K$ of
edges in the 1-skeleton of this complex.
While there are alternate formats in which
a knot and a 3-manifold may be presented,
all reasonable ones appear to be transformable to one another in
polynomial time.

In time polynomial in
$t$ we can check that the link of each vertex is connected and has Euler characteristic two,
which means that it is a sphere, and that the link of each edge
is a connected curve. These are the necessary conditions
to ensure that the underlying space of the complex is a 3-manifold $M$.
Similarly we can check in time polynomial in $t$ that the edges
of $K$ form a simple closed curve in $M$,
and that this curve represents a trivial element of the first homology
group of $M$ with integer coefficients.
We then form the second barycentric subdivision of the triangulation
of $M$, replacing each tetrahedron by 576 tetrahedra.
Removing all closed tetrahedra that meet $K$
results in a 3-manifold $M_K$ with a single torus boundary component,
the ``peripheral torus'' that surrounds the knot $K$.
The knot $K$ bounds a surface of genus $g$ in $M$
if and only if there is a surface
in $M_K$ of genus $g$ with a single
boundary component that is an essential curve on $\partial M_K$.
We restate our problem in the triangulated manifold $M_K$ as the
question of whether there exists in $M_K$ an orientable surface of genus
$g$ with a single essential boundary component on  $\partial M_K$.
By Schubert \cite{Sch:61}, if such a surface exists, then such a surface exists
among the fundamental normal surfaces in $M_K$.
The certificate consists of an integer
vector $\bw$ in ${\ZZ^{7t}}$
giving the normal surface coordinates of this surface.

Recall that not all vectors in ${\ZZ^{7t}}$ correspond to normal surfaces.
It is necessary that $\bw$ satisfies the matching, positivity and
quadrilateral conditions.
These conditions can be checked
in time which is linear in $t $.
So we can verify that there is a normal surface $F$
with $\bw = \bv (F)$.

To verify that $\partial F$ is essential on the peripheral torus,
we include in the certificate a non-trivial cycle in the
1-skeleton of $\partial M_K$ that intersects $\partial F$ in an odd number
of points. Such a cycle can be found in the 1-skeleton
of $\partial M_K$, since curves embedded in
the 1-skeleton generate its first homology.
Odd intersections with such a cycle
implies that $\partial F$ is non-separating  on $\partial M_K$,
and in particular does not bound a disk on $\partial M_K$,
ensuring that $\partial F$ corresponds to a longitude curve parallel to
$K$ in $N_K$.
Using the orbit counting algorithm, we can count the
number of components of the normal surface $F$ and verify that
$F$ is connected in time polynomial in $t \log W$,
where $W$ is the weight of $F$ the number of points in which $F$
meets the 1-skeleton of $M_K$. 
We apply Corollary~\ref{normalcount} to verify that
$\partial F$ is connected in time bounded by a polynomial in $t \log W$.
(This last step can be avoided. The number of essential boundary components
must be odd, since $\partial F$ represents a non-trivial element
in $H_1(\partial M_K; Z_2)$, and an even number of them
can be removed by joining adjacent pairs of curves with annuli on
$\partial M_K$.  This gives a surface with one boundary component and
the same Euler characteristic. Inessential boundary curves can be capped
while increasing Euler characteristic, which gives lower genus.)

To check that $F$ is orientable, we take the vector $2\bv (F) $
that doubles each coordinate of the normal surface $F$,
and apply the orbit counting algorithm to determine if
the corresponding normal surface  $\tilde F$ is connected.
Since $M$ is orientable, $\tilde F$ is connected
if $F$ is connected and non-orientable, and has
two components if $F$ is connected and orientable. 
Thus we can verify if $F$ is orientable in time which is
linear in $t \log W$.

The Euler Characteristic  $\chi (F0$
is determined by the number of vertices, edges and faces of $F$, which
are computable from its normal coordinates in time which is linear in $t$.
Following \cite{JT:95},
we let $t_i$ be the number of tetrahedra containing edge $e_i$
and set $\epsilon_{ij} = 1$ if the edge $e_i$ meets the $j^{th}$ normal disk.
Then the normal surface $F$ with coordinates $\{v_j \}$
has $\chi (F) = (1/2) f_3 - \sigma(F) + wt(F)$ where 
$\sigma(F) = \sum v_j$ and $wt(F) = \sum_{i,j} \epsilon_{ij} x_i/t_i $.
The values of $\epsilon_{ij}, t_i$ are determined by the triangulation
and are independent of $F$.  They can be computed in time polynomial in $t$.
Since $F$ is a connected orientable 
surface with one boundary component the
genus of $F$ is $(1 - \chi) /2$.  Thus we can determine the
genus of $F$ in time which is linear in $t $.

Theorem~\ref{bound} implies that the normal coordinates
of this surface are at most $t2^{7t+2}$.
There are $7t$ normal coordinates, and
each represents a triangle or quadrilateral,
so that the total number of intersections with the
1-skeleton satisfies $N \le 28{t^2}2^{7t+2}$.
In particular, $\log N$ is bounded above by a polynomial in $t$.
So the fact that $F$ is a spanning surface for $K$ can
be verified in time polynomial in $t$.
\qed

\section{An extended counting algorithm} \label{extended}

In this section we develop a generalized version
of the orbit counting algorithm, that counts not only
the number of orbits of a collection of isometries
between subintervals of an interval,
but also more general quantities which are useful in applications.
For example, we can use the extended algorithm to 
answer the following question: 
Given a normal surface and a triangulation,
how many times does each component of the surface 
intersect a fixed edged of the triangulation?
The extended algorithm allows one to effectively compute the normal
coordinates and the
Euler characteristic of each
connected component of the surface, hence the genus,  even when there are exponentially
many components. 
To carry out such computations, we extend the previous analysis to
pairings of weighted intervals, in which each point of the interval
has associated to it a vector in $\ZZ^d $. We are interested
in the sum of these vectors over an orbit, the {orbit weights}.

Consider again a pseudogroup of interval isometries acting on $[1,\dots,N]$.
We will assume that there is given as input
a nonnegative weight function $z: [1,\dots,N] \to \ZZ_+d $,
associating to each element of $[1,\dots,N]$
a vector in $\ZZ_+^d$ satisfying the following condition:
the weight at successive points $j$ and $j+1$ changes at most $4k$ times.
In our application, $k$ will be the initial number of pairings, and
$4k$ gives an upper bound on how many times an endpoint of a domain or
range of a pairing is reached as one moves across  $[1,\dots,N]$.
The algorithm proceeds as before, while maintaining 
data on the orbit weights.
We keep track of the orbit weights by maintaining two lists of weighted
subintervals.  The first
$$
L = \{([p_1,q_1],z_1),\dots,([p_m,q_m],z_m)\},\ 
1 \le p_i \le q_i < p_{i+1} 
$$
records the current weight values at each point in $[1,N]$,
with $q_m$ initially equal to $N$. This list
is updated as the algorithm proceeds.
Points in  the interval $[p_j,q_{j}]$ have constant
weight $z_j \in \ZZ^d$, and there are at most $4k$ such intervals. 
A second list of $t$ subintervals
$$
L' = \{([r_1,s_1],v_1),\dots,([r_t,s_t],v_t)\} , \ 1 \le r_i \le s_i < r_{i+1}
$$
consists of a collection of intervals $[ r_i , s_i]$
paired with a vector $v_i \in \ZZ_+^d$.
This pair represents $(s_i - r_i +1 ) $ orbits, one for each point in the interval
$[ r_i , s_i]$, and to each of these orbits is assigned the orbit weight $v_i$.
Initially empty, at the algorithm's conclusion
$L'$ records the total number of orbits
$ \displaystyle \sum_i (s_i - r_i +1 ) $,
along with the orbit weight $v_i$ assigned to
each of the $(s_i - r_i +1 ) $ orbits in the interval $[r_i, s_i]$.

We define an additional operation, called {\em transferring
weights  by a pairing $g$}.  Suppose that $g : [a, b]\to [c , d]$
is a pairing and that
$[c , d]$ carries $n$ different weights, given by the list
$$
\{([c = r_1,s_1],v_1),\dots,([r_n,s_n=d],v_n)\} , \ 1 \le r_i \le s_i < r_{i+1} .
$$
The weight function can be split into
constant functions on $n$ subintervals of $[c , d]$, where $1 \le n \le d-c+1$.

The transfer operation sets the weights on $[c , d]$ to zero and
keeps the orbit weights the same by translating
the weight vectors of $[c , d]$
to smaller orbit representatives, as below:

Case 1: $g$ is orientation preserving and $b < c$.
Set the weights on $[c , d]$ to zero and
for each $1 \le j \le n$, add $v_j$ to $g^{-1}([r_j,s_j])$.

Case 2:
$g$ is orientation preserving and $b \ge  c$.
Then $g : [a, b] \to [c, d]$ is a periodic pairing of period $t=c-a$.
We set the weights of points in $[c, d]$ to zero and adjust weights
in $[a, c-1]$ to preserve the orbit weight.
The points in $[c, d]$ have weights given by the intersection
of $[c, d]$ with intervals in $L$.
These weights are described by
$$
\{([c = r_1,s_1],v_1),\dots,([r_n,s_n=d],v_n)\} , \ 1 \le r_i \le s_i < r_{i+1} .
$$
For each interval $([r_j,s_j],v_j),\ 1 \le j \le n$, of width $w_j$ and constant weight $v_j$,
add $\displaystyle [[\frac{w_j}{t}]]v_j $
to the weight of each point in $[a,c-1]$, and add an additional $v_j$ to
the weight of each point in $[a,c-1]$ that is
congruent mod$(t)$ to a point in $\displaystyle [a +  [[\frac{w_j}{t}]]t,s_j]$, if any
such point exists.

Case 3: 
$g$ is orientation reversing.
We first trim $g$. This does not affect the orbits or the orbit weights.
For a trimmed, orientation reversing pairing 
$g : [a, b] \to [c, d]$,
set the weights on $[c , d]$ to zero and
for each $1 \le j \le n$ add $v_j$ to $g^{-1}([r_j,s_j])$.

\begin{lemma} \label{orbitweights.same}
The operation of transfer 
sets the weights on $[c, d]$ to zero and 
preserves the orbit weights of a collection of pairings.
The number of distinct weights taken by the weight function on $[1,N]$ increases
by at most four following a transfer operation.
\end{lemma}
\bproof
In Cases (1) and (3), the decrease in the weight function
at one point $x$ in an orbit is
exactly offset by an equal increase at the point
$g^{-1}(x)$ in the same orbit.

In Case (2), each point in $[a, d]$ is
in the orbit of a unique point in $[a,t-1] = [a, c-1]$ under the iterates of $g$.
The total weight of an orbit
within a  periodic pairing is transferred to the orbit representative
in this initial subinterval by applying powers of $g$.
Adding weight $v_j$ to an orbit representative in $[a,c-1]$ of each point
in $[r_j,s_j$ while setting the weight at that point to zero, preserves the
orbit weight. The resulting weights on $[a,c-1]$ are gotten by adding
appropriate multiples of $v_j$, with the factor being the number of orbits of
a point in $[a,c-1]$ that lie in $[r_j,s_j]$.   
The number of orbits under $g$ of a point in
$[a, c-1]$ that lie in $[r_j,s_j]$ is
$ \displaystyle [[\frac{w_j}{t}]]$ 
or one more than this for points whose orbit hits the last $w_j  (\mbox{mod } t)$ points
of $[r_j,s_j]$.
It follows that the transfer operation in Case (2) preserves the number of orbits and the orbit weight.

Setting the weights on $[c,d]$
to zero can cause at most two new points where a weight change occurs.
The transferred weights from a constant weight interval $[r_j, s_j]$
result in a net increase of
at most two pairs of successive points where the weight changes, at the
preimages of its two endpoints.
In Cases (1) and (3), transferred weights from
an interval with non-constant weights results
in an increase in the number of weight changes 
in the domain of $g$.  The increase in the number of weight changes in the 
domain is exactly canceled by the decrease in the number of weight changes in
the range of $g$, except possibly for two extra weight changes at the boundary points
$a, b$ of the domain.
In Case (2) the same holds, but with $[a,c-1]$ replacing the domain.
In each case the number of constant weight intervals $m'$ is increased by at most
four during a transfer operation.
\qed

We now describe the modified algorithm.  Again $N'$ represents
the current interval length, and we set $m'$ to be the
current number of constant weight intervals. \\

\noindent
{\bf Weighted orbit counting algorithm:}
Let $\{ g_i : [a_i, b_i] \to [c_i, d_i],\  1 \le i \le k \}$
be a collection of pairings between subintervals of $\{ 1, 2, \dots N \}$,
and let
$$
L = \{([p_1,q_1],z_1),\dots,([p_m,q_m],z_m)\} ,\  1 \le p_i \le q_i < p_{i+1} \le  N
$$
be a list representing a collection of weights on $[1,N]$, with the weight 
on $[p_j, q_j]$ equal to $z_j \in \ZZ^d $.
Initialize a second weight list $L'$ to be empty.
The algorithm proceeds as before, reducing the interval size
$N$ until it reaches zero, but this time
keeping track of orbit weights by maintaining the lists $L,L'$.
\begin{enumerate}
\item 
Search through the pairings and delete any pairings which are restrictions
of the identity.
Leave the weight lists $L,L'$ unchanged.

\item 
Search for and contract static intervals.
If the interval $[r,s]$ is contracted, has constant weight $z$,
and is not contained in a
larger contracted interval with the same weight $z$, add an interval
$([N'+1, N'+s-r+1], z)$ of
width $s-r+1$ to the end of $L'$, with associated weight $z$.
Alter $L$ by replacing $[1,N']$ by $[1, N'-(s-r+1)]$,
and altering each $g_i$ by replacing any point $x$ in the
domain or range of $g_i$ with $x>s$ by $x-(s-r+1)$.
Replace the weight function
by a new weight function $w'$,
which at points $x>s$
satisfies $w'(x) = w(x+(s-r+1))$ and agrees with $w$ at points $x<r$.

\item 
Trim all orientation reversing pairings whose domain and range
overlap.  Leave the weight lists $L,L'$ unchanged.

\item 
Search for pairs of periodic pairings 
$g_i$ and $g_j$ whose domains and ranges
satisfy the condition of  Lemma~\ref{merge2}.
If any such pair exists
then perform a merger as in Lemma~\ref{merge2},
replacing $g_i$ and $g_j$ by a single periodic pairing,
with translation distance $\mathrm{GCD}(t_i, t_j)$,
acting on the union of the domains and ranges of $g_i$ and $g_j$.
Leave the weight lists $L,L'$ unchanged.
Repeat until not mergers can be performed.

\item 
Find a maximal $g_i$. For each $g_j$ with $j \ne i$, if the range
of $g_j$ is contained in $[c_i, N']$, transmit $g_j$ by $g_i$.
Leave the weight lists $L,L'$ unchanged.

\item 
Find the smallest value of $c$ such that the interval $[c, N']$
intersects the range of at most one pairing $g_i$, with
$g_i : [a_i, b_i] \to [c_i, N']$.
Transfer the weights on $[c_i, N']$ by $g_i$, and then
truncate the pairing $g_i$.

\item 
If the interval size $N'$ has decreased to zero, output the list $L'$
and stop.
Otherwise start again with Step (1).
\end{enumerate}

For a  $\ZZ^d$-valued function on $[1,2, \dots, N]$ whose values are 
given by the list
$ L = \{([p_1,q_1],z_1),\dots,([p_m,q_m],z_m)\} ,\  1 \le p_i \le q_i < p_{i+1}$,
define the {\em total weight}
of $L$ to be $\sum_{i=1}^m |z_i|$.

\begin{theorem} \label{weightedcount}
Suppose there is a pseudogroup generated by $k$ pairings with $\ZZ^d$-valued weights,
$\{g_i\}_{i=1}^k$ on $[1,N]$ such that there is a partition of $[1,N]$ into
$m$ disjoint subintervals in which the weights are constant, and such that the total
weight is at most $D$.
Then the weighted orbit counting algorithm outputs a list with one point for
each orbit and corresponding orbit weights,
and runs in time polynomial in $ km d \log D \log N$.
\end{theorem}
\bproof
We will check that the running time of the algorithm
is larger than that of the previous unweighted version
by a factor which is a polynomial in $ md \log D $.

The proof that the algorithm terminates is the same as that
given for Theorem~\ref{orbitcount}. There is some extra overhead
involved in keeping track of weights that modifies the calculation
of the running time. We indicate these additional calculations below.
We now check that at each step in the algorithm the orbit weight is unchanged for
any orbit remaining in $L$, and that eliminated orbits have their
orbit weights correctly recorded in $L'$.

As we run through the steps of the algorithm,
Steps (1),(3),(4),(5) and (7) preserve the orbit structure, the
weight function and the interval $[1,N']$,
so neither of the lists $L,L'$ are changed.
The number of constant weight sub-intervals is also unchanged.

Step (2), contraction, does change the orbit structure, and also
shortens $[1,N']$. In Step (2) the procedure adds
the eliminated orbits and their weights to $L'$.
The number of constant weight intervals is not increased,
and may be decreased.
Maintaining the two lists requires at most $O(md \log D) $ additional
steps.

In Step (6), truncation, points are eliminated from the end of the
interval $[1,N']$.  However since we first transfer the weights
of these points, the eliminated points all have weight
zero and the weight of an orbit is unaffected.

The number of steps involved in resetting the weights in $L$
for a transfer operation is given by a polynomial in $m' d \log D \log N$.
Since $m'$ increases by at most four at each of the polynomially many steps of the algorithm,
$m'$ is bounded by a polynomial in $m k \log N$.

Combining the running time of each of the steps, whose number is given
by a polynomial in $ k \log N$, gives a 
polynomial in $ km d \log D \log N$ for the total running time.
\qed

\begin{corollary} \label{eulercount}
Let $M$ be a 3-manifold with a triangulation $T$
containing $t$ tetrahedra and
let $F$ be a normal surface in $M$ 
of total weight $W$.
There is a procedure for counting the number of components of $F$
and determining the topology of each component
which runs in time polynomial in $t \log W$.
\end{corollary}
\bproof
We begin as in Corollary~\ref{normalcount}
by assigning an integer in $[1,N]$ to each point of intersection
between the normal surface and an edge of the triangulation,
where $N$ is the total number of intersections of $F$ with the 1-skeleton,
and again associate three pairings to each face of the triangulation,
one to each pair of edges in the face.
The number of pairings that results is bounded above by $12t$.

We next define a weight function $w(x)$ which assigns
integer weights $(z_1, z_2, \dots, z_{7t})$
to each point in $[1,N]$.
Initially $z_i$ is set to zero for all $i$ at all points $x \in [1,N]$.
A tetrahedron can have as many as five distinct elementary disk types
with non-zero coefficients, four triangles and one quadrilateral.
If the $j^{th}$ elementary disk type occurs, then
fix one of the edges that it meets, and add $1$ to the $j^{th}$
component of the weight vector at
each of the indices that the $j^{th}$ elementary disk meets on that edge.
The orbit weights are then the normal
coordinates of the components of the normal surface $F$.
Each point in the output list $L'$ corresponds to a component of the
normal surface with normal coordinates given by the corresponding
weight in $\ZZ^{7t}$.
Theorem~\ref{weightedcount}
tells us that the list $L'$ is computed in time
polynomial in $ km d \log D \log N$.
We now bound these constants in terms of $t$.

Since each edge of a tetrahedron meets at most three disk types
in that tetrahedron, each edge of a tetrahedron can contribute at
most six points at which the weight vector changes.  Given six edges to
each tetrahedron, we have $m \le 36t$.
As before we have a bound for the number of pairings $k \le 12t$ and the number
of normal coordinates is given by $d = 7t$.
The total weight bounds the normal coordinates, and with $D=W$ we get a bound on 
the running time that is polynomial in $t \log W$.
\qed

\begin{corollary}
Let $M$ be a 3-manifold with a triangulation $T$
containing $t$ tetrahedra and
let $F$ be a fundamental normal surface in $M$.
There is a procedure for counting the number of components of $F$
and the topology of each component
which runs in time polynomial in $t$.
\end{corollary}
\bproof
Theorem~\ref{bound} gives a bound for the normal coordinates of $F$ of
$t2^{7t+2}$. 
Recall that there are $7t$ normal coordinates, and
each represents a triangle or quadrilateral,
so the total number of intersections with the
1-skeleton satisfies $W \le 28{t^2}2^{7t+2}$.
In particular, $\log W$ is bounded above by a polynomial in $t$.
Plugging this in for $W$ in Corollary~\ref{eulercount} we get a bound
for the running time which is a polynomial in $t$.
\qed

\section{The complexity of minimal spanning area} \label{area}

In this section we examine the complexity of
the problem of determining the  smallest area of a
spanning surface for a curve in a 3-dimensional manifold.  Such
an area calculation problem seems at first to be ill suited to
a complexity analysis, since it has real solutions
depending on a choice of Riemannian metric.

We recast the area calculation problem into
a discretized form where its complexity can be analyzed.
Given a curve in a suitably discretized Riemannian 3-manifold,
we ask whether it bounds a
surface of area less than $C$, where $C$ is an integer.
To describe a metric on a 3-manifold with a finite amount
of data, we restrict to piecewise flat metrics,
and manifolds constructed from collections of flat tetrahedra and triangular
prisms whose faces are identified by isometries. The curvature
of such PL metrics can be defined as a limit of smooth curvatures,
and is concentrated along their edges and vertices.
A particular manifold in this class is described by a decomposition into
tetrahedra or triangular prisms with
a rational (or integer) length assigned to each edge. The metric on this
tetrahedron or prism is then taken as the metric on the Euclidean
tetrahedron or prism with those edge lengths. In the case of a prism
we also set the angles of quadrilateral faces to be right angles.
Prisms are allowed in this construction in order to form metrics with rational lengths
on spaces that are products.
Identified 2-dimensional faces are required to be isometric.
We do not require that
the total angle around an edge is $2\pi$,
nor do we make any metric conditions at a vertex.
This type of metric is described by a
finite set of data, and can be used to approximate
Riemannian metrics on a manifold.
Up to scaling, we can take all the edge lengths to be integers.
We call these objects {\em metrized PL 3-manifolds}.
A curve is given as a collection of edges in the
1-skeleton of $M$.
We will show that given an integer $C$, determining whether the smallest spanning surface
for a curve in such a 3-manifold has area less than $C$
is {\bf NP}-hard.

\begin{tabbing}
{\em Problem:} MINIMAL-SPANNING-AREA\\
{\em INSTANCE:} A 3-dimensional metrized PL manifold $M$, a 1-dimensional \\
curve $K$ in the 1-skeleton of $M$, and a natural number $C$. \\ 
{\em QUESTION:} \=Does the curve bound a
surface of area $A \le C$?
\end{tabbing}

The size of an instance is given by the number of bits needed to describe all the edge lengths
and $C$.

\begin{theorem} \label{nphard}
MINIMAL-SPANNING-AREA is {\bf NP}-hard.
\end{theorem}
 
\bproof
We reduce in polynomial time
an instance of the  {\bf NP}-hard problem ONE-IN-THREE SAT
to an instance of MINIMAL-SPANNING-AREA. This shows that
MINIMAL-SPANNING-AREA is at least as hard, up to polynomial time reduction,
as  ONE-IN-THREE SAT.

As a first step, we set up a 2-dimensional
version of MINIMAL-SPANNING-AREA.
We then construct a 3-manifold by a thickening process,
with the property that a minimizing surface must remain within the 
2-complex.

\begin{figure}[hbtp]
\centering
\includegraphics[width=.8\textwidth]{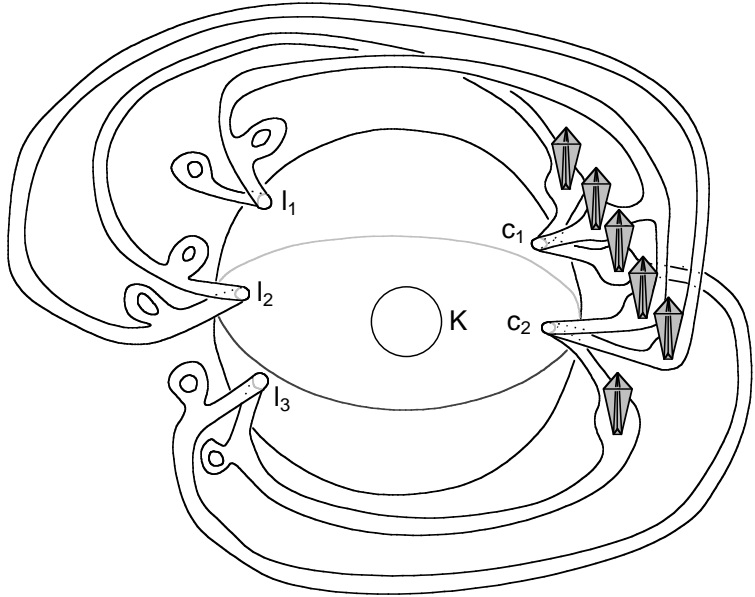}
\caption{A metrized branching surface corresponding to the
boolean expression
$(x_1 \vee x_2 \vee x_3) \wedge (x_1 \vee \bar x_2 \vee \bar x_3)$.
The picture is not to scale.
The shaded prisms are constructed to each have area one,
while the rest of the surface has total area less than 1/2.}
\label{metrized}
\end{figure}

Given a boolean expression representing
an instance of ONE-IN-THREE SAT, we construct
a triangulated metrized 2-complex and an integer $C$.
This complex contains a curve $K$ 
with the property that the expression admits a satisfying
assignment if and only if $K$ bounds a surface of area less than $C$.
This metrized complex is shown in Figure~\ref{metrized}
for the expression $(x_1 \vee x_2 \vee x_3) \wedge (x_1 \vee \bar x_2 \vee \bar x_3)$.

The branching surface is similar to the one used in the proof of
Theorem~\ref{thmnph}, but carries the additional structure of a
metrized triangulation, whose triangles have
flat metrics of prescribed edge length.
The metrized triangles are constructed so that near each
of the $m$ boundary components corresponding to clauses of
the boolean expression there are three triangulated disks of area close to one,
one on each of the three handles coming into the punctured sphere
near the boundary component.
These disks are shaded in Figure~\ref{metrized}.
Each of these shaded disks is chosen
to have area between 1 and $1 + 1/2m$.
The surface is constructed so that
the union of all triangles in the rest of the surface has total area less than 1/2.

We saw in Theorem~\ref{thmnph} that a
spanning surface which has minimal genus
goes over each of the shaded disks at most once and goes
over exactly one shaded disk for each of the $m$ clauses.
It follows that such a surface has total area $m <  A < m+1$.
Furthermore, a satisfying
assignment for ONE-IN-THREE SAT leads to an embedded
spanning surface, with the satisfying values of the variables
selecting branches of the surface, and such a
spanning surface has area less than $m+1$.
So an instance
of ONE-IN-THREE SAT can be reduced to an instance of 
MINIMAL-SPANNING-AREA for this 2-complex.

To pass to a 3-manifold, thicken each triangle
in the branching surface to a triangular prism, triangulated
as in the proof of Theorem~\ref{thmnph}, and with 
a product metric. This produces a 3-manifold $M$ which is
a thickened up version of the 2-complex.
Projection to the branching surface is area non-increasing, and area
decreasing for a surface with boundary on the branching surface
but not contained in it. Therefore a least
area surface spanning $K$ must lie on the branching surface.
A closed manifold $DM$ with a piecewise-smooth metric
can be obtained by a doubling construction
as in Theorem~\ref{thmnph}. The doubling involution is an isometry, so that
reflecting a surface meeting $DM \backslash M$ into $M$ does not increase area.
It follows that the embedded spanning surface on the branching surface
is a least area surface in $DM$.
\qed

\section{Open questions}

Among many unresolved questions are:

1. Does determining knot genus remain {\bf NP}-hard
if we restrict to knots in the 3-sphere? 

2. Is determining the genus of a knot in a 3-manifold {\bf NP}?
This amounts to showing that finding a lower bound to the knot's genus
is an {\bf NP} problem, in contrast to the upper bound we have investigated.
Recall that the genus of a knot is the least possible
genus of all spanning surfaces.
We have shown that certifying that the genus is at most $g$ is {\bf NP},
but left open the possibility that the genus may be smaller than $g$.
If the answer to this question is yes,
then we can certify that a non-trivial knot has
positive genus, and it would follow that
{\sc UNKNOTTING} is both {\bf NP} and {\bf coNP}.

\end{document}